\documentclass[12pt]{article}
\usepackage{graphicx}
\usepackage{amsmath,amsthm,amssymb,enumerate}
\usepackage{euscript,mathrsfs}
\usepackage{color}
\usepackage{dsfont}
\usepackage[left=2cm,right=2cm,top=3.5cm,bottom=3.5cm]{geometry}
\usepackage{color}
\usepackage[framemethod=tikz]{mdframed}
\allowdisplaybreaks

\usepackage{soul}

\catcode`\@=11 \@addtoreset{equation}{section}

\catcode`\@=12

\newtheorem{Theorem}{Theorem}[section]
\newtheorem{Proposition}[Theorem]{Proposition}
\newtheorem{Lemma}[Theorem]{Lemma}
\newtheorem{Corollary}[Theorem]{Corollary}

\theoremstyle{definition}
\newtheorem{Definition}[Theorem]{Definition}

\newtheorem{Remark}[Theorem]{Remark}

\newcommand{\bTheorem}[1]{
	\begin{Theorem} \label{T#1} }
	\newcommand{\eT}{\end{Theorem}}

\newcommand{\bProposition}[1]{
	\begin{Proposition} \label{P#1}}
	\newcommand{\eP}{\end{Proposition}}

\newcommand{\bLemma}[1]{
	\begin{Lemma} \label{L#1} }
	\newcommand{\eL}{\end{Lemma}}

\newcommand{\bCorollary}[1]{
	\begin{Corollary} \label{C#1} }
	\newcommand{\eC}{\end{Corollary}}

\newcommand{\bRemark}[1]{
	\begin{Remark} \label{R#1} }
	\newcommand{\eR}{\end{Remark}}

\newcommand{\bDefinition}[1]{
	\begin{Definition} \label{D#1} }
	\newcommand{\eD}{\end{Definition}}

\newcommand{\phiB}{\bfphi_{\mathcal{B}}}

\newcommand{\Ds}{\mathbb{D}_x}

\newcommand{\tvm}{\widetilde{\vc{m}}}

\newcommand{\tS}{\widetilde{S}}
\newcommand{\bfphi}{\boldsymbol{\varphi}}

\newcommand{\bFormula}[1]{
	\begin{equation} \label{#1}}
	\newcommand{\eF}{\end{equation}}

\newcommand{\Ov}[1]{\overline{#1}}

\newcommand{\vr}{\varrho}

\newcommand{\tvr}{\widetilde{\vr}}
\newcommand{\tvu}{{\tilde \vu}}
\newcommand{\tvt}{\tilde \vt}

\newcommand{\vt}{\vartheta}
\newcommand{\vu}{\vc{u}}
\newcommand{\vm}{\vc{m}}

\newcommand{\vc}[1]{{\bf #1}}

\newcommand{\Div}{{\rm div}_x}
\newcommand{\Grad}{\nabla_x}

\newcommand{\dx}{\,{\rm d} {x}}

\newcommand{\dt}{\,{\rm d} t }

\newcommand{\intO}[1]{\int_{\Omega} #1 \ \dx}

\newcommand{\D}{{\rm d}}

\newcommand{\ep}{\varepsilon}

\newcommand{\vtB}{\vt_B}

\newcommand{\br}{ \nonumber \\ }

\def\softd{{\leavevmode\setbox1=\hbox{d}%
		\hbox to 1.05\wd1{d\kern-0.4ex{\char039}\hss}}}
\definecolor{Cgrey}{rgb}{0.85,0.85,0.85}
\definecolor{Cblue}{rgb}{0.50,0.85,0.85}
\definecolor{Cred}{rgb}{1,0,0}
\definecolor{fancy}{rgb}{0.10,0.85,0.10}

\newcommand\Cbox[2]{%
	\newbox\contentbox%
	\newbox\bkgdbox%
	\setbox\contentbox\hbox to \hsize{%
		\vtop{
			\kern\columnsep
			\hbox to \hsize{%
				\kern\columnsep%
				\advance\hsize by -2\columnsep%
				\setlength{\textwidth}{\hsize}%
				\vbox{
					\parskip=\baselineskip
					\parindent=0bp
					#2
				}%
				\kern\columnsep%
			}%
			\kern\columnsep%
		}%
	}%
	\setbox\bkgdbox\vbox{
		\color{#1}
		\hrule width  \wd\contentbox %
		height \ht\contentbox %
		depth  \dp\contentbox
		\color{black}
	}%
	\wd\bkgdbox=0bp%
	\vbox{\hbox to \hsize{\box\bkgdbox\box\contentbox}}%
	\vskip\baselineskip%
}

\mdfdefinestyle{MyFrame}{%
	linecolor=black,
	outerlinewidth=1pt,
	roundcorner=5pt,
	innertopmargin=\baselineskip,
	innerbottommargin=\baselineskip,
	innerrightmargin=10pt,
	innerleftmargin=10pt,
	backgroundcolor=white!20!white}

\begin{document}


\title{Unconditional stability of equilibria in thermally driven compressible fluids}

\author{Eduard Feireisl
		\thanks{The work of E.F. was supported by the
			Czech Sciences Foundation (GA\v CR), Grant Agreement
			24--11034S. The Institute of Mathematics of the Academy of Sciences of
			the Czech Republic is supported by RVO:67985840. } \and Yong Lu \thanks{The work of Y. L.  was partially supported by the  NSF of China under Grant 12171235.} \and Yongzhong Sun \thanks{The work of Y. S. was partially supported by the NSF of China under Grant 12071211.}}

\date{}

\maketitle

\medskip

\centerline{Institute of Mathematics of the Academy of Sciences of the Czech Republic}

\centerline{\v Zitn\' a 25, CZ-115 67 Praha 1, Czech Republic. Email: feireisl@math.cas.cz.}

\medskip

\centerline{Department of Mathematics, Nanjing University}

\centerline{22 Hankou Road, 210093 Nanjing, China. Email: luyong@nju.edu.cn.}

\medskip

\centerline{Department of Mathematics, Nanjing University}

\centerline{22 Hankou Road, 210093 Nanjing, China.  Email: sunyz@nju.edu.cn.}

\begin{abstract}
	
	We show that small perturbations of the spatially homogeneous equilibrium of a thermally driven compressible viscous 	fluid are globally stable. Specifically, any weak solution of the evolutionary Navier--Stokes--Fourier system driven by thermal convection converges to an equilibrium as time goes to infinity. The main difficulty to overcome is the fact the problem does not admit any obvious Lyapunov function. The result applies, in particular, to the 	Rayleigh--B\' enard convection problem. 
	
\end{abstract}




\section{Introduction}
\label{p}

Understanding the long time behaviour of thermally driven fluid flows is crucial in the mathematical theory of turbulence, see Constantin et al. \cite{Const1}, 
\cite{CoFoMaTe},  Davidson \cite{DAVI}, Foias et al. \cite{FMRT} among others.  An iconic example is the Rayleigh--B\' enard problem, where 
the fluid is confined between two parallel plates heated from below and subjected to the gravitational force, see Bormann \cite{Borm2}, \cite{Borm1}, Cao et al. 
\cite{CaJoTuWh}. The onset of turbulence is usually preceded by the loss of stability of the equilibrium solution triggered by reaching critical values 
of the Reynolds number and other parameters as the case may be. 

We consider the motion of a viscous compressible and heat conducting fluid driven by an inhomogeneous temperature 
field prescribed on the boundary of the physical domain $\Omega \subset R^3$.  
The time evolution of the mass density 
$\vr = \vr(t,x)$, the temperature $\vt = \vt(t,x)$, and the velocity $\vu = \vu(t,x)$ of the fluid is described by the 
\emph{Navier--Stokes--Fourier (NSF) system} of field equations:
\begin{align} 
	\partial_t \vr + \Div (\vr\vu) &= 0, \label{p1} \\
	\partial_t (\vr \vu) + \Div (\vr \vu \otimes \vu) + \Grad p(\vr,\vt) &= \Div \mathbb{S}(\vt, \Ds \vu) + \vr \Grad G, \label{p2} \\ 
	\partial_t (\vr e(\vr, \vt)) + \Div (\vr e(\vr, \vt) \vu) + \Div \vc{q}(\vt, \Grad \vt) &= \mathbb{S}(\vt, \Ds \vu): \Ds \vu - p(\vr, \vt) \Div \vu,
	\label{p3}
	\end{align}	
where $	\Ds \vu = \frac{1}{2} \left( \Grad \vu + \Grad^t \vu \right)$. The pressure $p = p(\vr, \vt)$ is related to the internal energy $e = e(\vr, \vt)$ through
\emph{Gibbs' equation}
\begin{equation} \label{p4}
	\vt D s = D e + p D \left( \frac{1}{\vr} \right),
	\end{equation}
where $s = s(\vr, \vt)$ is the entropy. 
The viscous stress tensor $\mathbb{S} = \mathbb{S}(\vt, \Ds \vu)$ is given by
\emph{Newton's rheological law} 
\begin{equation} \label{p5}
	\mathbb{S} (\vt, \Ds \vu) = \mu(\vt) \left( \Grad \vu + \Grad \vu^t - \frac{2}{3} \Div \vu \mathbb{I} \right) + 
	\eta(\vt) \Div \vu \mathbb{I}. 
	\end{equation}
Finally, the heat flux $\vc{q} = \vc{q}(\vt, \Grad \vt)$ is given by 
\emph{Fourier's law} 
\begin{equation} \label{p6}
	\vc{q} = - \kappa(\vt) \Grad \vt.
	\end{equation}

We suppose the fluid occupies a bounded domain $\Omega \subset R^3$ and impose the standard no--slip boundary conditions for the velocity, 
\begin{equation} \label{p7}
	\vu|_{\partial \Omega} = 0.
	\end{equation}
Our main goal is to study the long time behaviour of the system driven by the boundary temperature distribution 
\begin{equation} \label{p8}
	\vt|_{\partial \Omega} = \vtB. 
	\end{equation}
	
If $G = 0$ and $\vt_B = \Ov{\vt}_B$ is a positive constant, the NSF system \eqref{p1}--\eqref{p3} admits a unique (constant) 
\emph{static} solution 
$\vr_s = m_0 |\Omega|^{-1}$, $\vt_s = \Ov{\vt}_B$, $\vu_s = 0$, where
$$
m_0 = \intO{ \vr }
$$
is the total mass of the fluid.
We consider a small perturbation of the boundary temperature 
\begin{equation} \label{p11a}
	\| \vtB - \Ov{\vt}_B \|_{C^{2 + \nu}(\partial \Omega)} \leq \ep 
\end{equation} 
together with a small potential driving force 
\begin{equation} \label{p12}
	\| G \|_{C^1(\Ov{\Omega})} \leq  \ep,
\end{equation}
for some $\ep$ sufficiently small.  Under these conditions, the NSF system admits a \emph{stationary} state $(\vr_s, \vt_s, \vu_s)$
solving the 
\emph{stationary problem}
\begin{align}
	\Div (\vr_s \vu_s) &= 0, \label{p9}\\
	\Div (\vr_s \vu_s \otimes \vu_s) + \Grad p(\vr_s, \vt_s) &= \Div \mathbb{S}(\vt_s, \Ds \vu_s) + \vr_s \Grad G,\ \vu_s|_{\partial \Omega} = 0, \label{p10} \\
	\Div (\vr_s e(\vr_s, \vt_s) \vu_s ) - \Div (\kappa (\vt_s) \Grad \vt_s ) &= \mathbb{S}(\vt_s, \Ds \vu) : \Ds \vu - 
	p(\vr_s, \vt_s) \Div \vu_s ,\ \vt_s|_{\partial \Omega} = \vt_B, \label{p11}\\  
	\mbox{with the prescribed total mass}\ &\intO{\vr_s} = m_0 > 0  .  \label{p11b}
\end{align}
Moreover, the stationary solution $(\vr_s, \vt_s, \vu_s)$ is a small perturbation of the static equilibrium:
\begin{align}
	\left\| \vr_s - \frac{m_0}{|\Omega|} \right\|_{C^1(\Ov{\Omega})} \lesssim \ep ,  \quad \left\| \vt_s - \Ov{\vt}_B 	\right\|_{C^1(\Ov{\Omega})} \lesssim \ep , \quad \left\| \vu_s \right\|_{C^2(\Ov{\Omega}; R^3)} \lesssim \ep. 
	\label{p14}
\end{align}
Here and hereafter, the symbol $a \lesssim b$ mean there exists a positive constant $C$ such that 
$a \leq C b$. We shall also write $a \approx b$ whenever $a \lesssim b$ and $b \lesssim a$. Note carefully that for a non--constant boundary 
temperature $\vtB$, the stationary velocity $\vu_s$ need not vanish, cf. e.g. Daniels et al. \cite{DBPB}.

Under assumptions \eqref{p11a} and \eqref{p12}, the existence of small stationary solutions in the class \eqref{p14} was established e.g. by Valli and Zajaczkowski \cite[Section 5]{VAZA}. In particular, Valli and Zajaczkowski \cite[Sections 4, 5]{VAZA} show  that small and smooth perturbations of the equilibrium are \emph{locally} stable in the class of smooth solutions.
Our goal is to extend this result to \emph{global stability} in the class of weak solutions. Specifically, we show 
\[
(\vr, \vt, \vu) (t, \cdot) \to (\vr_s, \vt_s, \vu_s) \ \mbox{as}\ t \to \infty 
\]  
for \emph{any} global--in--time weak solution of the NSF system. The unconditional stability of equilibrium is rather surprising as the problem does not admit any explicit Lyapunov function.  To the best of our knowledge, this is the first result concerning \emph{global} stability of equilibrium  of the NSF system driven by thermal convection.

We point out that the Navier--Stokes--Fourier system with the boundary conditions \eqref{p7}, \eqref{p8} is an open 
(non--conservative) dissipative system, where convergence to equilibrium is not obvious and, in general, possibly not true. This is in sharp contrast with the closed system endowed with the \emph{conservative} boundary conditions
\[
\vu|_{\partial \Omega} = 0,\quad \Grad \vt \cdot \vc{n}|_{\partial \Omega} = 0,
\]
for which the unconditional convergence to a unique equilibrium was established in \cite{FP20}, see also the monograph 
\cite{FeiPr}.

\medskip

The main result of the present paper is stated in Theorem \ref{TM1} in Section \ref{M} below.
The following are the main ingredients of the proof:
\begin{itemize}
	\item We use the framework of \emph{weak solutions} introduced recently in \cite{ChauFei}, \cite{FeiNovOpen}, 
	see Section  \ref{w}. This concept is based, in particular, on incorporating the entropy inequality 
	as a part of the definition in accordance with the dissipative character of the fluid flow. 
	
	\item As shown in \cite{FeiSwie}, the Navier--Stokes--Fourier system endowed with the boundary conditions \eqref{p7}, \eqref{p8} is dissipative in the sense of Levinson, meaning it admits a bounded absorbing set, see Section \ref{A}. 
	In particular, the energy of any global--in--time solution is bounded by a universal constant. 
	
	\item The relative energy introduced in \cite{ChauFei} comparing the distance between a weak solution and the 
	equilibrium solution represents a suitable candidate for a Lyapunov function. Unfortunately, in the present setting, 
	we recover only the relative energy inequality containing non--trivial forcing term, see Section \ref{rie}. 
	
	\item The bulk of the paper is the analysis of the interaction between the forcing and the dissipative terms 
	appearing in the relative energy inequality, see Section \ref{ef}. The main difficulty is the absence of the 
	density in the dissipation term that must be compensated by a careful analysis of the pressure in the momentum 
	equation. This being established, the proof of convergence can be completed by standard arguments. 
	
	\item The result applies in a straightforward manner to the Rayleigh--B\' enard convection problem considered 
	in the class of compressible fluids, see Section \ref{cr}. 
	
	\end{itemize}

\section{Weak solutions to the Navier--Stokes--Fourier system}
\label{w}

We recall the concept of \emph{weak solution} to the NSF system introduced in \cite{ChauFei}, \cite{FeiNovOpen}. 

\subsection{Structural hypotheses concerning constitutive relations}
\label{S}

The mathematical theory of the NSF system developed in \cite{FeiNovOpen} requires certain restrictions to be imposed on the constitutive 
relations. These are mostly physically grounded hypotheses and we refer the reader to \cite[Chapter 1]{FeNo6A} for their physical background.

\subsubsection{Equation of state}

The equations of state relating the pressure and the internal energy to the state variables $(\vr, \vt)$ are that of a monoatomic gas augmented by the effect of thermal radiation:
\begin{align} 
	p(\vr, \vt) &= \vt^{\frac{5}{2}} P \left( \frac{\vr}{\vt^{\frac{3}{2}}} \right) + \frac{a}{3} \vt^4, \ a > 0, \label{S1} \\
	e(\vr, \vt) &= \frac{3}{2} \frac{ \vt^{\frac{5}{2}} }{\vr} P \left( \frac{\vr}{\vt^{\frac{3}{2}}} \right) + \frac{a}{\vr} \vt^4, \label{S2} 
	\end{align}
and, in accordance with Gibbs' relation \eqref{p4}, 
\begin{equation}
	s(\vr, \vt) = \mathcal{S} \left( \frac{\vr}{\vt^{\frac{3}{2}}} \right) + \frac{4a}{3 \vr} \vt^3, \label{S3}
	\end{equation}
where
\begin{equation}  
	\mathcal{S}'(Z) = -\frac{3}{2} \frac{ \frac{5}{3} P(Z) - P'(Z) Z }{Z^2}. \label{S4}
\end{equation}

Next, we impose the \emph{hypothesis of thermodynamics stability}
\[
\frac{\partial p(\vr, \vt)}{\partial \vr} > 0,\ \frac{\partial e(\vr, \vt)}{\partial \vt} > 0,
\]
which, expressed in terms of $P\in C^1[0,\infty)$, reads 
\begin{equation}
		 P(0) = 0,\ P'(Z) > 0 \ \mbox{for all}\ Z \geq 0,\ 0<\frac{ \frac{5}{3} P(Z) - P'(Z) Z}{Z} \leq c, \ \mbox{for all}\ Z > 0 \label{S5}
	\end{equation}
for some positive constant $c$.  In addition, there are extra hypotheses on the behaviour of $P$ in the degenerate area $Z = \frac{\vr}{\vt^{\frac{3}{2}}} \gg 1$. 
First, we impose the Third law of thermodynamics requiring
\begin{equation} 	
	\lim_{Z \to \infty} \mathcal{S}(Z) = 0, \label{S7}
\end{equation}
cf. Belgiorno \cite{BEL1}, \cite{BEL2}.
Second, we suppose 
	\begin{equation}\label{convex-pm}
	P \in C^1[0,\infty) \ \mbox{is convex}.
\end{equation}

\begin{Remark} \label{rP1}
 Note that $P(Z) = {\rm R}Z$ for moderate values of the degeneracy parameter $Z$ corresponds to the standard Boyle-Mariotte equation of state for gases. The convex deviation postulated in \eqref{convex-pm} is characteristic for Fermi gas formed by free electrons in the degenerate area $Z\gg 1$, see M\" uller and Ruggeri \cite[Chapter 4]{MURU}. As observed in \cite{Fei23a}, convexity of $P$ plays a crucial role in the proof of positivity of the temperature in the NSF system.
	
	\end{Remark}

Now, it follows from \eqref{S5} that the function 
\[
Z \mapsto \frac{ P(Z) }{Z^{\frac{5}{3}}} 
\]
is decreasing, and, consequently,
\begin{equation}
		\lim_{Z \to \infty} \frac{ P(Z) }{Z^{\frac{5}{3}}} = p_\infty \geq 0. 
		\label{S6a}
		\end{equation} 
Next, we observe that convexity of $P$ postulated in \eqref{convex-pm} implies $p_\infty > 0$. Indeed 
the Third law enforced through \eqref{S7}, together with \eqref{S4}, imply
\begin{equation} \label{S66}
\mathcal{S}(Z) = \frac{3}{2} \int_Z^\infty \frac{\frac{5}{3} P(s) - P'(s) s }{s^2} \ \D s \to 0 \ \mbox{as} \ Z \to \infty. 
\nonumber\end{equation}
A direct computation yields
\[
\int_Z^\infty \frac{\frac{5}{3} P(s) - P'(s) s }{s^2} \ \D s = - \int_Z^\infty s^{\frac{2}{3}}  \frac{\D }{\D s} \frac{P(s)}{s^{\frac{5}{3}}} \ \D s
\geq - Z^{\frac{2}{3}} \int_Z^\infty   \frac{\D }{\D s} \frac{P(s)}{s^{\frac{5}{3}}} \ \D s.
\]
Thus if $p_\infty = 0$ in \eqref{S6a} we conclude 
\[
\mathcal{S}(Z) \geq \frac{3}{2} Z^{\frac{2}{3}} \frac{P(Z)}{Z^{\frac{5}{3}}} = \frac{3}{2} \frac{P(Z)}{Z} \to 0 \ \mbox{as}\ Z \to \infty,
\]
which is impossible as $P'(0) > 0$ and $P$ is convex. We may therefore infer that 
\begin{equation}
	\lim_{Z \to \infty} \frac{ P(Z) }{Z^{\frac{5}{3}}} = p_\infty > 0.
	 \label{S6}
\end{equation}

Write 
\begin{equation} \label{form-p}
P(Z) = p_\infty Z^{\frac{5}{3}} + P_m(Z). 
\end{equation}
In accordance with \eqref{S4}, 
\[
S'(Z) = - \frac{3}{2} \frac{ \frac{5}{3} P_m(Z) - P'_m(Z) Z }{Z^2} < 0;
\]
whence, exactly as above and in view of \eqref{S6}, we get 
\begin{equation} \label{S91}
\frac{ P_m(Z) }{Z^{\frac{5}{3}}} \searrow 0 \ \mbox{as}\ Z \to \infty.
\nonumber \end{equation}
In particular, as $P_m(0) = 0$, $P'_m(0) = P'(0) > 0$ we get $P_m(Z) > 0$ for all $Z > 0$.
Finally, repeating the arguments leading to \eqref{S6}, we get 
\begin{equation} \label{S92}
	\mathcal{S}(Z) \geq \frac{3}{2} \frac{P_m(Z)}{Z} \to 0 \ \mbox{as}\ Z \to \infty.
\end{equation}
and, in view of \eqref{S5}, 
\begin{equation} \label{S93}	
\limsup_{Z \to \infty} P'_m(Z) = 0.
\end{equation}
Finally, we exploit the upper bound in hypothesis \eqref{S5}. Differentiating \eqref{form-p} and using \eqref{S6} we get 
\begin{equation} \label{S94}
P'_m(Z) = P'(Z) - \frac{5}{3} p_\infty Z^{\frac{2}{3}} \geq P'(Z) - \frac{5}{3} \frac{P(Z)}{Z^{\frac{5}{3}}} Z^{\frac{2}{3}} = 
\frac{P'(Z)Z - \frac{5}{3} P(Z) }{Z} \geq -c.
\end{equation}
It follows from \eqref{form-p}--\eqref{S94} that the pressure can be written in the form 
\begin{equation} \label{S9}
	p(\vr, \vt) = p_\infty \vr^{\frac{5}{3}} + \vt^{\frac{5}{2}} P_m \left( \frac{\vr}{\vt^{\frac{3}{2}}} \right) + \frac{a}{3} \vt^4, \quad a>0, \ p_\infty>0, 
\end{equation}
where the function $P_m \in C^1[0,\infty)$ enjoys the properties  
\begin{equation} \label{pm-infty-2}
P_m \geq 0, \ \lim_{Z \to \infty} \frac{P_m(Z)}{Z} = 0, 
\end{equation}
\begin{equation} \label{pm-infty-3}
- c\leq \liminf_{Z \to \infty}	P'_m(Z) \leq \limsup_{Z \to \infty} P'_m (Z) = 0. 
\end{equation}

\subsubsection{Transport coefficients}

As for the transport coefficients $\mu$, $\eta$, and $\kappa$, we suppose they are continuously differentiable functions 
of the temperature $\vt \in [0, \infty)$ satisfying
\begin{align}
	0 < \underline{\mu} \left(1 + \vt \right) &\leq \mu(\vt) \leq \Ov{\mu} \left( 1 + \vt \right),\
	|\mu'(\vt)| \lesssim 1, \label{S10} \\ 
	0 &\leq  \eta(\vt) \leq \Ov{\eta} \left( 1 + \vt \right)
	\label{S11}, \\ 
	0 < \underline{\kappa} \left(1 + \vt^\beta \right) &\leq  \kappa(\vt) \leq \Ov{\kappa} \left( 1 + \vt^\beta \right) \ \mbox{for some}\  \beta > 6.	
	\label{S12}
\end{align}

\subsection{Weak solutions}

We introduce the concept of global--in--time weak solutions to the NSF system. As our main goal is to study their 
long--time behavior, we shall completely ignore the initial data and consider the weak solutions defined on an open time interval $(T, \infty)$.

	\begin{Definition} \label{Dw1} {\bf (Global in time weak solutions)}
	
	Let $\Omega \subset R^3$ be a bounded Lipschitz domain. 
	A quantity $(\vr, \vt, \vu)$ is  \emph{weak solution} of the NSF system
	\eqref{p1}--\eqref{p6}, with the boundary conditions \eqref{p7}, \eqref{p8} in the 
	time--space cylinder $(T, \infty) \times \Omega$, $T \geq - \infty$,
	if the following holds:
	\begin{itemize}
		\item{\bf Equation of continuity.}
		$\vr \in L^\infty_{\rm loc}(T, \infty; L^{\frac{5}{3}} (\Omega))$, $\vr \geq 0$, and the integral
		identity
		\begin{equation} \label{w1}
			\int_T^\infty \intO{ \Big( \vr \partial_t \varphi + \vr \vu \cdot \Grad \varphi \Big) } \dt = 0
		\end{equation}
		holds for any $\varphi \in C^1_c((T, \infty) \times \Ov{\Omega})$. 	
		In addition, the renormalized version of \eqref{w1}
		\begin{equation} \label{w2}
			\int_T^\infty \intO{ \left( b(\vr) \partial_t \varphi + b(\vr) \vu \cdot \Grad \varphi +
				\Big( b(\vr) - b'(\vr) \vr \Big) \Div \vu \varphi \right) } \dt = 0
		\end{equation}
		holds for any $\varphi \in C^1_c((T, \infty) \times \Ov{\Omega})$ and any $b \in C^1(R)$, $b' \in C_c(R)$.
		
		\item{\bf Momentum equation.}
		
		$\vr \vu \in L^\infty_{\rm loc}(T, \infty ; L^{\frac{5}{4}}(\Omega; R^3))$,
		$\vu \in L^2_{\rm loc}(T, \infty; W^{1,2}(\Omega; R^3))$, $\vu|_{\partial \Omega} = 0$,
			and the integral identity
		\begin{align}
			\int_T^\infty &\intO{ \Big( \vr \vu \cdot \partial_t \bfphi + \vr \vu \otimes \vu : \Grad \bfphi + 
				p(\vr, \vt) \Div \bfphi \Big) } \br &= \int_T^\infty \intO{ \mathbb{S}(\vt, \Grad \vu) : \Grad \bfphi } \dt-
			\int_T^\infty \intO{ \vr \Grad G \cdot \bfphi } \dt  	
			\label{w3}
		\end{align}
		for any $\bfphi \in C^1_c((T, \infty) \times {\Omega}; R^3)$.

		\item {\bf Entropy inequality.}
		$\vt \in L^\infty_{\rm loc}(T, \infty; L^4(\Omega)) \cap L^2_{\rm loc}(T, \infty; 
		W^{1,2}(\Omega))$, $\vt > 0$ a.e. in $(T, \infty) \times \Omega$,
		$\log (\vt) \in L^2_{\rm loc}(T, \infty; W^{1,2}(\Omega))$, 
		\[
		\vt|_{\partial \Omega} = \vtB.
		\]
		The integral inequality
		\begin{align}
			\int_T^\infty &\intO{ \left( \vr s (\vr, \vt) \partial_t \varphi + \vr s (\vr, \vt) \vu \cdot \Grad \varphi + \frac{\vc{q} (\vt, \Grad \vt) }{\vt} \cdot
				\Grad \varphi \right) } \dt \br & \leq  - \int_T^\infty \intO{ \frac{\varphi}{\vt}
				\left( \mathbb{S}(\vt, \Grad \vu) : \Grad \vu - \frac{\vc{q} (\vt, \Grad \vt) \cdot \Grad \vt}{\vt} \right) }
			\dt	
		\nonumber
		\end{align}
		holds for any $\varphi \in C^1_c((T, \infty) \times {\Omega})$, $\varphi \geq 0$.
		
		\item {\bf Ballistic energy inequality.}
		The inequality
		\begin{align}
			&\int_T^\infty \partial_t \psi	\intO{ \left( \frac{1}{2} \vr |\vu|^2 + \vr e(\vr, \vt)  - \tvt \vr s(\vr, \vt) \right) } \dt \br
			\quad&
			- \int_T^\infty \psi \intO{ 	\frac{\tvt}{\vt} \left( \mathbb{S}(\vt, \Grad \vu) : \Grad \vu - \frac{\vc{q}(\vt, \Grad \vt) \cdot \Grad \vt}{\vt}  \right) } \dt \br
			& \geq  \int_T^\infty \psi \intO{\left( \vr s (\vr, \vt) \partial_t \tvt + \vr s (\vr, \vt) \vu \cdot \Grad \tvt + \frac{\vc{q}(\vt, \Grad \vt)}{ \vt} \cdot \Grad \tvt                    \right)      } \dt\br	
			&\quad -  \int_T^\infty \psi \intO{ \vr \Grad G \cdot \vu } \dt 
			\nonumber 
		\end{align}
		holds for any $\psi \in C^1_c(T, \infty)$, $\psi \geq 0$, and any
		$\tvt \in C^1((T, \infty) \times \Ov{\Omega})$, $\tvt > 0$, $\tvt|_{\partial \Omega} = \vtB$.

	\end{itemize}		
	
\end{Definition}

Under the hypotheses \eqref{S1}--\eqref{S7}, and \eqref{S10}--\eqref{S12}, 
the existence of global--in--time weak solutions for any finite energy initial data was proved in 
\cite{ChauFei} (see also \cite{FeiNovOpen} for a more detailed proof).

\section{Dissipativity in the sense of Levinson}
\label{A}

Here and hereafter, we suppose the boundary temperature as well as the gravitational potential are independent of time: 
\begin{equation} \label{A1}
	\vtB = \vtB(x), \ G = G(x). 
	\end{equation}
Recall that 
\begin{equation} \label{A2}
	m_0 = \intO{ \vr } 
	\end{equation}
denotes the total mass of the fluid, which is a constant of motion even in the class of weak solutions. 

Next, we introduce a quantity measuring the amplitude of the data, 
\begin{align} 
	\| (\mbox{data}) \| &= m_0 + m_0^{-1} + \|G \|_{W^{1,\infty}( \Omega)} +  
	\| \vtB^{-1} \|_{L^\infty(\partial \Omega)}  + \| \vtB  \|_{W^{2, \infty}(\partial \Omega)},
	\label{A3}
\end{align}
together with the total energy $E$, 
\[
E(\vr, \vu, \vt) = \frac{1}{2} \vr |\vu|^2 + \vr e(\vr, \vt).
\]

The following result was proved in \cite[Theorem 3.1]{FeiSwie}, see also \cite[Theorem 3.1]{BreFei23}. 

\begin{Proposition}[\bf Bounded absorbing set] \label{PA1}
	
	Let $\Omega \subset R^3$ be a bounded domain of class $C^{2 + \nu}, \ \nu > 0$. Suppose the thermodynamic 
	functions $p$, $e$, and $s$ satisfy the hypotheses \eqref{S1}--\eqref{convex-pm}, while the transport coefficients 
	$\mu$, $\eta$, and $\kappa$ comply with \eqref{S10}--\eqref{S12}.

		Then there exists a positive quantity $\mathcal{E}_\infty = \mathcal{E}_\infty	( \| (\mbox{\rm data}) \| ) $ -- a bounded 
	non--decreasing function of the norm of the data specified in \eqref{A3} -- such that 
	\begin{equation} \label{A5}
		\limsup_{t \to \infty} \intO{ E(\vr, \vt, \vu)(t, \cdot) } \leq \mathcal{E}_\infty 
		\end{equation} 
for any weak solution of the NSF system defined on $(T, \infty) \times \Omega$ in the sense of Definition \ref{Dw1}.	
	\end{Proposition}

In view of the specific form of the equation of state, notably \eqref{S1}, \eqref{S6}, we deduce from Proposition \ref{PA1}
that 
\begin{equation} \label{A6} 
{\rm ess} \limsup_{t \to \infty} \left( \| (\vr \vu) (t, \cdot) \|_{L^{\frac{5}{4}}(\Omega; R^3)} + 
\| \vr(t, \cdot) \|_{L^{\frac{5}{3}}(\Omega)} + \| \vt (t, \cdot) \|_{L^4(\Omega)} \right) \leq \mathcal{F}_\infty,	
	\end{equation}
where $\mathcal{F}_\infty$ is a universal constant depending solely on  $\| (\mbox{data}) \|$.

\section{Main result -- unconditional stability}
\label{M}

We are ready to state our main result concerning unconditional stability of small perturbations of the constant static solutions. A suitable ``distance functional'' is provided by the \emph{relative energy} 
\begin{align}
	E &\left( \vr, \vt, \vu \Big| \tvr , \tvt, \tvu \right) \br &= \frac{1}{2}\vr |\vu - \tvu|^2 + \vr e - \tvt \Big(\vr s - \tvr s(\tvr, \tvt) \Big)- 
	\Big( e(\tvr, \tvt) - \tvt s(\tvr, \tvt) + \frac{p(\tvr, \tvt)}{\tvr} \Big)
	(\vr - \tvr) - \tvr e (\tvr, \tvt) \br &= 
	\frac{1}{2}\vr |\vu - \tvu|^2 + \vr e - \tvt \vr s - \Big( e(\tvr, \tvt) - \tvt s(\tvr, \tvt) + \frac{p(\tvr, \tvt)}{\tvr} \Big) \vr + p(\tvr, \tvt).
\nonumber
\end{align}
It follows from the hypothesis of thermodynamics stability that
the relative energy interpreted in terms of the conservative entropy variables $(\vr, S = \vr s, \vm = \vr \vu)$ represents a Bregman distance associated to the convex energy functional 
\[
E(\vr, S, \vm) = \frac{1}{2} \frac{|\vm|^2}{\vr} + \vr e(\vr, S).
\]
see  \cite{FeiNov20}, \cite[Chapter 3, Section 3.1]{FeiNovOpen}. Specifically, 
\[
E \left( \vr, S, \vm \Big| \tvr , \tS, \tvm \right) = E(\vr, S, \vm) - \left< \partial E(\tvr, \tS, \tvm) ; (\vr - \tvr, S - \tS, \vm - \tvm) \right> - E(\tvr, \tS, \tvm).
\]

\begin{Theorem}[\bf Unconditional stability] \label{TM1}
Let $\Omega \subset R^3$ be a bounded domain of class $C^{2 + \nu},\ \nu > 0$. Suppose the thermodynamic 
functions $p$, $e$, and $s$ satisfy the hypotheses  \eqref{S1}--\eqref{convex-pm}, while the transport coefficients 
$\mu$, $\eta$, and $\kappa$ comply with \eqref{S10}--\eqref{S12}. Let the data $G, \vtB, m_0$ satisfy
\begin{equation} \label{M3}
	\| G \|_{C^1(\Omega)} \leq \ep, \quad \left\| \vt_B - \Ov{\vt}_B \right\|_{C^{2+\nu}(\Ov{\Omega})}\leq \ep,\quad  0 < m_0 < \infty 
	\end{equation}	
for some constant $\Ov{\vt}_B >0$. 

Then there exists $\ep_0 > 0$ such that
\begin{equation} \label{M4}
	\intO{ E \left( \vr, \vt, \vu \Big| \vr_s , \vt_s, \vu_s \right) (t, \cdot) } \to 0 
\ \mbox{as}\ t \to \infty
\end{equation}
whenever $0< \ep \leq \ep_0$ for any global--in--time weak solution $(\vr, \vt, \vu)$ of the NSF system defined on $(T, \infty)$,
with the total mass
\[
\intO{\vr } = m_0,
\]
where $(\vr_s, \vt_s, \vu_s)$ is the solution of the stationary problem \eqref{p9}--\eqref{p11b}.

	\end{Theorem}
	
\begin{Remark} \label{RR1}
	The convergence in terms of the relative energy stated in \eqref{M4} along with the uniform bounds \eqref{A6} yields, in particular, 
	\begin{align}
	\vr(t, \cdot) &\to \vr_s \ \mbox{as}\ t \to \infty \  \mbox{in}\ L^{q_1}(\Omega),\ 1 \leq q_1 < \frac{5}{3}, \br 
	\vr \vu (t, \cdot) &\to \vr_s \vu_s \ \mbox{as}\ t \to \infty \  \mbox{in}\ L^{q_2}(\Omega; R^3),\ 1 \leq q_2 < \frac{5}{4}, \br
	 \vr s(\vr, \vt) (t, \cdot) &\to \vr_s s(\vr_s, \vt_s) \ \mbox{as}\ t \to \infty \  \mbox{in}\ L^{q_3}(\Omega),\ 1 \leq q_3 < \frac{4}{3}.
		\nonumber
	\end{align}
	
	\end{Remark}	

The rest of the paper is devoted to the proof of Theorem \ref{TM1}. To begin, in view of the result of Valli and Zajackowski \cite{VAZA}, 
we observe that the stationary solution $(\vr_s, \vt_s, \vu_s)$ is smooth and small satisfying \eqref{p14}.

\section{Proof of Theorem \ref{TM1}}
\label{c}

Let $(\vr_s, \vt_s, \vu_s)$ be the stationary solution to \eqref{p9}--\eqref{p11b}  satisfying \eqref{p14}.  Of course, the solution $(\vr_s, \vt_s,\vu_s)$ depends on $\ep$ but we omit this to keep the notation simple. Note that 
existence and uniqueness of the stationary solution in the class \eqref{p14} was established by Valli and Zajaczkowski \cite{VAZA}.	

\subsection{Relative entropy inequality}
\label{rie}

We recall the relative energy inequality, see \cite[Chapter 12, Section 12.3.2]{FeiNovOpen}:
\begin{align}
	&\frac{\D}{\dt} \intO{ E \left( \vr, \vt, \vu \ \Big| \tvr , \tvt, \tvu \right) } 
	\br&	+  \intO{ 	\frac{\tvt}{\vt} \left( \mathbb{S}(\vt, \Grad \vu) : \Grad \vu - \frac{\vc{q}
			(\vt, \Grad \vt) \cdot \Grad \vt}{\vt} \right) }  \br
	&\quad \leq -  \intO{ \Big( \vr (\vu - \tvu) \otimes (\vu - \tvu)   + p (\vr, \vt) \mathbb{I} - \mathbb{S}(\vt, \Grad \vu) \Big) : \Ds \tvu  }  \br
	&\quad \quad + \intO{ \frac{\vr}{\tvr} (\vu - \tvu) \cdot \Grad p(\tvr, \tvt) }  \br
	&\quad \quad-  \intO{ \vr \Big(  \partial_t \tvu + \tvu \cdot \Grad \tvu + \frac{1}{\tvr} \Grad p(\tvr, \tvt) -
		\Grad G \Big) \cdot (\vu - \tvu)     }  \br
	&\quad \quad -  \intO{\left( \vr \Big( s(\vr, \vt)  - s(\tvr, \tvt) \Big) \partial_t \tvt + \vr \Big(s (\vr, \vt)  - s(\tvr, \tvt) \Big) \vu \cdot \Grad \tvt + \frac{\vc{q} (\vt, \Grad \vt) }{ \vt} \cdot \Grad \tvt                    \right)      } \br
	&\quad \quad +  \intO{ \left( \left(1 - \frac{\vr}{\tvr} \right) \partial_t p(\tvr, \tvt) - \frac{\vr}{\tvr} \vu \cdot \Grad p(\tvr,\tvt)      \right)     } 
	\ \mbox{in}\ \mathcal{D}'(T, \infty),
	\label{c2}
\end{align}
where $(\vr, \vt, \vu)$ is a weak solution to the NSF system and $(\tvr, \tvt, \tvu)$ is an arbitrary trio of sufficiently smooth test functions satisfying the boundary conditions \eqref{p7}, \eqref{p8}. The natural idea is to use the the stationary equation $(\vr_s, \vt_s, \vu_s)$ as test functions and examine the resulting inequality. 
Thus, we plug $(\tvr, \tvt, \tvu) = (\vr_s, \vt_s, \vu_s)$ in \eqref{c2} obtaining
\begin{align}
	&\frac{\D}{\dt} \intO{ E \left( \vr, \vt, \vu \ \Big| \vr_s , \vt_s, \vu_s \right) } 	+  \intO{ 	\frac{\vt_s}{\vt}  \mathbb{S}(\vt, \Grad \vu) : \Grad \vu -\mathbb{S}(\vt,\Grad \vu): \Grad\vu_s } \br 
	& +\intO{\left( \frac{\vc{q} (\vt, \Grad \vt) }{ \vt} \cdot \Grad \vt_s    - \frac{\vc{q}(\vt, \Grad \vt) \cdot \Grad \vt}{\vt} \right) }  \br 
	&\quad \leq  -  \intO{ \Big( \vr (\vu - \vu_s) \otimes (\vu - \vu_s)    \Big) : \Grad \vu_s  }  -\intO{ \left( p(\vr,\vt) - p(\vr_s,\vt_s)\right)\Div\vu_s} \br 
	&\quad\quad -\intO{ \frac{\vr}{\vr_s}\Div\mathbb{S}(\vt_s,\Grad\vu_s)\cdot (\vu-\vu_s) }-  \intO{\left(  \vr \Big(s (\vr, \vt)  - s(\vr_s, \vt_s) \Big) \vu \cdot \Grad \vt_s    \right)      } .
	\label{c4}
\end{align}
Furthermore, following the arguments of \cite[Chapter 12]{FeiNovOpen}, we rewrite \eqref{c4} in the form 
\begin{align}
	&\frac{\D}{\dt} \intO{ E \left( \vr, \vt, \vu \ \Big| \vr_s , \vt_s, \vu_s \right) }  +   \intO{  \mathbb{S}(\vt,\Grad (\vu -\vu_s))  : \Grad(\vu - \vu_s) }	\br
	&\qquad +  \intO{ \frac{\vt_s -\vt}{\vt}\mathbb{S}(\vt,\Grad\vu): \Grad\vu   +\frac{\vt - \vt_s}{\vt_s} \mathbb{S}(\vt_s, \Grad \vu_s) : \Grad \vu_s  } \br
	&\qquad + \intO{ \left[\frac {\vt_s - \vt    }{\vt} \frac{\kappa (\vt) |\Grad \vt|^2 }{\vt} + \frac {\vt - \vt_s }{\vt_s} \frac{\kappa (\vt_s) |\Grad \vt_s|^2 }{\vt_s}   + \left( \frac{\kappa (\vt) \Grad \vt }{\vt}  - \frac{\kappa (\vt_s) \Grad \vt_s }{\vt_s}      \right) \cdot
	\Big( \Grad \vt - \Grad \vt_s \Big) \right]} \br
	&\quad \leq  -  \intO{ \Big( \vr (\vu - \vu_s) \otimes (\vu - \vu_s)    \Big) : \Grad \vu_s  } - \intO{ \frac{\vr-\vr_s}{\vr_s}\Div\mathbb{S}(\vt_s,\Grad\vu_s)\cdot (\vu - \vu_s) } \br
	&\qquad - \intO{\left(  \vr \Big(s (\vr, \vt)  - s(\vr_s, \vt_s) \Big) (\vu-\vu_s) \cdot \Grad \vt_s    \right) - \left(  (\vr-\vr_s) \Big(s (\vr, \vt)  - s(\vr_s, \vt_s) \Big) \vu_s \cdot \Grad \vt_s    \right) } \br
	&\qquad + \intO{ \Div\vu_s \left( p(\vr_s,\vt_s) - \frac{\partial{p}(\vr_s,\vt_s)}{\partial\vr}(\vr_s - \vr) - \frac{\partial{p}(\vr_s,\vt_s)}{\partial\vt}(\vt_s - \vt) - p(\vr,\vt)\right) }\br
	&\qquad + \intO{ \left( s(\vr_s,\vt_s) - \frac{\partial{s}(\vr_s,\vt_s)}{\partial\vr}(\vr_s - \vr) - \frac{\partial{s}(\vr_s,\vt_s)}{\partial\vt}(\vt_s - \vt) - s(\vr,\vt)\right)\vu_s\cdot\Grad\vt_s }.
	\label{c5}
\end{align}

\subsection{Estimates involving temperature dissipation}\label{et}

For the dissipation term  related to the temperature in  \eqref{c5}, we deduce two relations: 
\begin{align}
	&\frac{\vt_s - \vt}{\vt}   \frac{\kappa (\vt) |\Grad \vt|^2 }{\vt} + \frac{\vt - \vt_s}{\vt_s} 
\frac{ \kappa (\vt_s) |\Grad \vt_s|^2 }{\vt_s} + 
\left(  \frac{\kappa (\vt) \Grad \vt }{\vt} - \frac{\kappa (\vt_s) \Grad \vt_s }{\vt_s}  \right) \cdot (\Grad \vt - \Grad \vt_s)  \br 
& = \vt_s \frac{\kappa(\vt)}{\vt^2}  |\Grad \vt - \Grad \vt_s|^2 + 
2 \frac{\kappa(\vt)}{\vt^2}  (\vt_s - \vt) (\Grad \vt - \Grad \vt_s) \cdot \Grad \vt_s \br 
&\quad + \left(  \frac{\kappa(\vt) }{\vt} - \frac{\kappa(\vt_s) }{\vt_s} \right) (\Grad \vt - \Grad \vt_s) \cdot \Grad \vt_s + \left( \frac{\kappa(\vt) }{\vt} - \frac{\kappa(\vt_s) }{\vt_s} \right) 
\left( \frac{\vt_s - \vt}{ \vt_s } \right) |\Grad \vt_s |^2  \br
&\quad + \frac{\kappa(\vt)}{\vt}(\vt_s - \vt) \left( \frac{1}{\vt} - \frac{1}{\vt_s} \right) |\Grad \vt_s |^2 \br
&\geq \vt_s \frac{\kappa(\vt)}{\vt^2}  |\Grad \vt - \Grad \vt_s|^2 + 
2 \frac{\kappa(\vt)}{\vt^2}  (\vt_s - \vt) (\Grad \vt - \Grad \vt_s) \cdot \Grad \vt_s \br 
&\quad + \left(  \frac{\kappa(\vt) }{\vt} - \frac{\kappa(\vt_s) }{\vt_s} \right) (\Grad \vt - \Grad \vt_s) \cdot \Grad \vt_s + \left( \frac{\kappa(\vt) }{\vt} - \frac{\kappa(\vt_s) }{\vt_s} \right) 
\left( \frac{\vt_s - \vt}{ \vt_s } \right) |\Grad \vt_s |^2,	
	\label{cc1}
\end{align}	
and
\begin{align} 
	&\frac{\vt_s - \vt}{\vt}   \frac{\kappa (\vt) |\Grad \vt|^2 }{\vt} + \frac{\vt - \vt_s}{\vt_s} 
\frac{ \kappa (\vt_s) |\Grad \vt_s|^2 }{\vt_s} + 
\left(  \frac{\kappa (\vt) \Grad \vt }{\vt} - \frac{\kappa (\vt_s) \Grad \vt_s }{\vt_s}  \right) \cdot (\Grad \vt - \Grad \vt_s)  \br 
&= \vt_s  \frac{\kappa (\vt)}{\vt^2}  |\Grad \vt|^2  +  \frac{\vt}{\vt_s} 
\frac{ \kappa (\vt_s) |\Grad \vt_s|^2 }{\vt_s} - \frac{\kappa (\vt) \Grad \vt }{\vt}  \cdot \Grad \vt_s - 
\frac{\kappa (\vt_s) \Grad \vt_s }{\vt_s} \cdot  \Grad \vt .
\label{cc2}
\end{align}

Next, we fix two constants $\underline{\vt}$ and $ \Ov{\vt}$ such that
\begin{equation} \label{cc3}
	0 < \underline{\vt} \leq \frac{1}{2} \inf_\Omega \vt_s \leq 2 \sup_\Omega \vt_s \leq \Ov{\vt}.
\end{equation}

\begin{enumerate}
	
	\item {\bf Estimates for intermediate $\vt$.}
	
	It follows from \eqref{cc1} and Cauchy-Schwartz inequality that
\begin{align}
\mathds{1}_{\underline{\vt} \leq \vt \leq \Ov{\vt} }	& \left[ \frac{\vt_s - \vt}{\vt}   \frac{\kappa (\vt) |\Grad \vt|^2 }{\vt} + \frac{\vt - \vt_s}{\vt_s} 
	\frac{ \kappa (\vt_s) |\Grad \vt_s|^2 }{\vt_s} + 
	\left(  \frac{\kappa (\vt) \Grad \vt }{\vt} - \frac{\kappa (\vt_s) \Grad \vt_s }{\vt_s}  \right) \cdot (\Grad \vt - \Grad \vt_s) \right]  \br 
	&\geq \frac{\vt_s}{2} \mathds{1}_{\underline{\vt} \leq \vt \leq \Ov{\vt} }\frac{\kappa(\vt)}{\vt^2}  |\Grad \vt - \Grad \vt_s|^2 -  \|\Grad \vt_s\|_{L^{\infty}(\Omega)}^2 c(\underline{\vt}, \Ov{\vt} ) \mathds{1}_{\underline{\vt} \leq \vt \leq \Ov{\vt} } |\vt - \vt_s |^2 ,	
	\label{cc4}
\end{align}
where $c(\underline{\vt}, \Ov{\vt} ) $ denotes a constant depending on $\underline{\vt}, \Ov{\vt} $.

\item {\bf Estimates for small $\vt$.}

It follows from \eqref{cc2} and algebraic H\" older inequality, 
\begin{align}
	\mathds{1}_{0 < \vt \leq \underline{\vt} }	& \left[ \frac{\vt_s - \vt}{\vt}   \frac{\kappa (\vt) |\Grad \vt|^2 }{\vt} + \frac{\vt - \vt_s}{\vt_s} 
	\frac{ \kappa (\vt_s) |\Grad \vt_s|^2 }{\vt_s} + 
	\left(  \frac{\kappa (\vt) \Grad \vt }{\vt} - \frac{\kappa (\vt_s) \Grad \vt_s }{\vt_s}  \right) \cdot (\Grad \vt - \Grad \vt_s) \right]  \br 
	&\geq \frac{2 \vt_s}{3} \mathds{1}_{0 < \vt \leq \underline{\vt} } \frac{\kappa(\vt)}{\vt^2}  |\Grad \vt |^2 -  \|\Grad \vt_s\|_{L^{\infty}(\Omega)}^2 c(\underline{\vt} ) \mathds{1}_{0 < \vt \leq \underline{\vt} } .	
	\label{cc5}
\end{align}
 where $c(\underline{\vt}) $ denotes a  constant depending on $\underline{\vt} $.
Using the fact
\begin{equation}\label{cc4-1}
 |\Grad \vt - \Grad \vt_s |^2 \leq 2 |\Grad \vt|^2 + 2 |\Grad \vt_s|^2,
\nonumber \end{equation}
and 
$$
\frac{\kappa(\vt)}{\vt} \geq \underline{\vt}^{-1} \underline{\kappa}  \quad \mbox{whenever $0<\vt<\underline{\vt}$},
$$
we derive from \eqref{cc5}  that
\begin{align}
	\mathds{1}_{0 < \vt \leq \underline{\vt} }	& \left[ \frac{\vt_s - \vt}{\vt}   \frac{\kappa (\vt) |\Grad \vt|^2 }{\vt} + \frac{\vt - \vt_s}{\vt_s} 
	\frac{ \kappa (\vt_s) |\Grad \vt_s|^2 }{\vt_s} + 
	\left(  \frac{\kappa (\vt) \Grad \vt }{\vt} - \frac{\kappa (\vt_s) \Grad \vt_s }{\vt_s}  \right) \cdot (\Grad \vt - \Grad \vt_s) \right]  \br 
	&\geq \frac{ \vt_s}{2} \mathds{1}_{0 < \vt \leq \underline{\vt} } \frac{\kappa(\vt)}{\vt^2}  |\Grad \vt |^2 
	+ \frac{ \vt_s}{2} \mathds{1}_{0 < \vt \leq \underline{\vt} } |\Grad \vt - \Grad \vt_s |^2 
	- \|\Grad \vt_s\|_{L^{\infty}(\Omega)}^2 c(\underline{\vt}, \Ov{\vt} ) \mathds{1}_{0 < \vt \leq  \underline{\vt} } .	
	\label{cc6}
\end{align}	
Finally, we realize the inequality
\begin{equation} \label{cc7}
\mathds{1}_{0 < \vt \leq  \underline{\vt} } \leq c(\underline{\vt} ) \mathds{1}_{0 < \vt \leq  \underline{\vt} } 
|\vt - \vt_s |^2.
\end{equation}	
Thus we may infer from \eqref{cc6} and \eqref{cc7} that 
\begin{align}
	\mathds{1}_{0 < \vt \leq \underline{\vt} }	& \left[ \frac{\vt_s - \vt}{\vt}   \frac{\kappa (\vt) |\Grad \vt|^2 }{\vt} + \frac{\vt - \vt_s}{\vt_s} 
	\frac{ \kappa (\vt_s) |\Grad \vt_s|^2 }{\vt_s} + 
	\left(  \frac{\kappa (\vt) \Grad \vt }{\vt} - \frac{\kappa (\vt_s) \Grad \vt_s }{\vt_s}  \right) \cdot (\Grad \vt - \Grad \vt_s) \right]  \br 
	&\geq \frac{ \vt_s}{2} \mathds{1}_{0 < \vt \leq \underline{\vt} } \frac{\kappa(\vt)}{\vt^2}  |\Grad \vt |^2 
	+ \frac{ \vt_s}{2} \mathds{1}_{0 < \vt \leq \underline{\vt} } |\Grad \vt - \Grad \vt_s |^2 
	- \|\Grad \vt_s\|_{L^{\infty}(\Omega)}^2 c(\underline{\vt} ) \mathds{1}_{0 < \vt \leq  \underline{\vt} } 
	|\vt - \vt_s |^2.	
	\label{cc8}
\end{align}	

\item {\bf Estimates for large $\vt$.}

Using \eqref{cc2} again we deduce:
\begin{align} 
\mathds{1}_{\Ov{\vt} \leq \vt  }	& \left[ \frac{\vt_s - \vt}{\vt}   \frac{\kappa (\vt) |\Grad \vt|^2 }{\vt} + \frac{\vt - \vt_s}{\vt_s} 
	\frac{ \kappa (\vt_s) |\Grad \vt_s|^2 }{\vt_s} + 
	\left(  \frac{\kappa (\vt) \Grad \vt }{\vt} - \frac{\kappa (\vt_s) \Grad \vt_s }{\vt_s}  \right) \cdot (\Grad \vt - \Grad \vt_s) \right] \br 
	&\geq \frac{3 \vt_s}{4} \mathds{1}_{\Ov{\vt} \leq \vt  }  \frac{\kappa (\vt)}{\vt^2}  |\Grad \vt|^2  
	- \|\Grad \vt_s\|_{L^{\infty}(\Omega)}^2 c(  \Ov{\vt} ) \mathds{1}_{\Ov{\vt} \leq \vt } 
	\big(1 + \kappa (\vt) \big). 
	\label{cc9}
\end{align}
 where $c(\Ov{\vt}) $ denotes a constant depending on $\Ov{\vt} $.
Next, similarly to the above,
\[
\frac{\kappa (\vt)}{\vt^2} |\Grad \vt - \Grad \vt_s |^2 \leq 2  \frac{\kappa (\vt)}{\vt^2} |\Grad \vt |^2 + 
2  \frac{\kappa (\vt)}{\vt^2} |\Grad \vt_s |^2,
\]
whence \eqref{cc9} reduces to 
\begin{align} 
	\mathds{1}_{\Ov{\vt} \leq \vt  }	& \left[ \frac{\vt_s - \vt}{\vt}   \frac{\kappa (\vt) |\Grad \vt|^2 }{\vt} + \frac{\vt - \vt_s}{\vt_s} 
	\frac{ \kappa (\vt_s) |\Grad \vt_s|^2 }{\vt_s} + 
	\left(  \frac{\kappa (\vt) \Grad \vt }{\vt} - \frac{\kappa (\vt_s) \Grad \vt_s }{\vt_s}  \right) \cdot (\Grad \vt - \Grad \vt_s) \right] \br 
	&\geq \frac{ 2\vt_s}{3} \mathds{1}_{\Ov{\vt} \leq \vt  }  \frac{\kappa (\vt)}{\vt^2}  |\Grad \vt|^2  
	+  \frac{ \vt_s}{2} \mathds{1}_{\Ov{\vt} \leq \vt  }  \frac{\kappa (\vt)}{\vt^2}  |\Grad \vt - \Grad \vt_s|^2
	- \|\Grad \vt_s\|_{L^{\infty}(\Omega)}^2 c( \Ov{\vt} ) \mathds{1}_{\Ov{\vt} \leq \vt } 
	 \kappa (\vt). 
	\label{cc10}
\end{align}

Finally, we write 
\begin{align}
\mathds{1}_{\Ov{\vt} \leq \vt } \kappa(\vt) &\approx  \mathds{1}_{\Ov{\vt} \leq \vt }  + \mathds{1}_{\Ov{\vt} \leq \vt } \vt^\beta  \approx   \mathds{1}_{\Ov{\vt} \leq \vt } +  \mathds{1}_{\Ov{\vt} \leq \vt } \left( \vt^{\frac{\beta}{2}} - \Ov{\vt}^{\frac{\beta}{2}} \right)^2 +   \overline\vt^\beta \br
 &\approx    \left(   \left[ \vt^{\frac{\beta}{2}} - \Ov{\vt}^{\frac{\beta}{2}} \right]^+ \right)^2 + \mathds{1}_{\Ov{\vt} \leq \vt } \left(1+\overline \vt^{\beta}\right) .
\label{cc11} 
\end{align}
With the choice of $\overline \vt$ in \eqref{cc3} we have
$$
\vt|_{\partial \Omega} = \vt_{B} \leq \sup_{\Omega} \vt_{B} \leq \sup_{\Omega} \vt_{s} \leq \frac{\overline\vt}{2}.
$$
This implies 
$$
\left[ \vt^{\frac{\beta}{2}} - \Ov{\vt}^{\frac{\beta}{2}} \right]^+ = 0 \quad \mbox{on} \ \partial \Omega.
$$
Using Poincar\' e inequality, we get 
\begin{equation} \label{cc12}
\intO{ \left(   \left[ \vt^{\frac{\beta}{2}} - \Ov{\vt}^{\frac{\beta}{2}} \right]^+ \right)^2 } 
\leq c_P \intO{ \mathds{1}_{\Ov{\vt} \leq \vt } |\Grad \vt^{\frac{\beta}{2}} |^2 } 
\lesssim \intO{ \mathds{1}_{\Ov{\vt} \leq \vt } \frac{\kappa (\vt)}{\vt^2} |\Grad \vt |^2 } ,
\end{equation}
where $c_{P}$ is the constant coming from the Poincar\'e inequality.  Recall $\|\Grad \vt_s\|_{L^{\infty}(\Omega)} \lesssim \ep  \leq  \ep_{0}$ with $\ep_{0}$ sufficiently small. 
Consequently, \eqref{cc10} together with \eqref{cc12} give rise to 
\begin{align} 
	\mathds{1}_{\Ov{\vt} \leq \vt  }	& \left[ \frac{\vt_s - \vt}{\vt}   \frac{\kappa (\vt) |\Grad \vt|^2 }{\vt} + \frac{\vt - \vt_s}{\vt_s} 
	\frac{ \kappa (\vt_s) |\Grad \vt_s|^2 }{\vt_s} + 
	\left(  \frac{\kappa (\vt) \Grad \vt }{\vt} - \frac{\kappa (\vt_s) \Grad \vt_s }{\vt_s}  \right) \cdot (\Grad \vt - \Grad \vt_s) \right] \br 
	&\geq \frac{ \vt_s}{2} \mathds{1}_{\Ov{\vt} \leq \vt  }  \frac{\kappa (\vt)}{\vt^2}  |\Grad \vt|^2  
	+  \frac{  \vt_s}{2} \mathds{1}_{\Ov{\vt} \leq \vt  }  \frac{\kappa (\vt)}{\vt^2}  |\Grad \vt - \Grad \vt_s|^2
	-\|\Grad \vt_s\|_{L^{\infty}(\Omega)}^2 c( \Ov{\vt} ) \mathds{1}_{\Ov{\vt} \leq \vt } . 
	\label{cc12-1}
\nonumber\end{align}
	
To conclude, we recall 
\begin{equation} \label{cc13}
\mathds{1}_{\Ov{\vt} \leq \vt } \leq c(\Ov{\vt}) \mathds{1}_{\Ov{\vt} \leq \vt  } 
|\vt - \vt_s |^2,
\nonumber \end{equation}
and therefore
\begin{align} 
	\mathds{1}_{\Ov{\vt} \leq \vt  }	& \left[ \frac{\vt_s - \vt}{\vt}   \frac{\kappa (\vt) |\Grad \vt|^2 }{\vt} + \frac{\vt - \vt_s}{\vt_s} 
	\frac{ \kappa (\vt_s) |\Grad \vt_s|^2 }{\vt_s} + 
	\left(  \frac{\kappa (\vt) \Grad \vt }{\vt} - \frac{\kappa (\vt_s) \Grad \vt_s }{\vt_s}  \right) \cdot (\Grad \vt - \Grad \vt_s) \right] \br 
	&\geq \frac{ \vt_s}{2} \mathds{1}_{\Ov{\vt} \leq \vt  }  \frac{\kappa (\vt)}{\vt^2}  |\Grad \vt|^2  
	+  \frac{  \vt_s}{2} \mathds{1}_{\Ov{\vt} \leq \vt  }  \frac{\kappa (\vt)}{\vt^2}  |\Grad \vt - \Grad \vt_s|^2
	- \|\Grad \vt_s\|_{L^{\infty}(\Omega)}^2 c( \Ov{\vt} ) \mathds{1}_{\Ov{\vt} \leq \vt } |\vt - \vt_s |^2 . 
	\label{cc13-1}
\end{align}

\item {\bf Conclusion.}

Putting the above estimates together, notably \eqref{cc4}, \eqref{cc8}, \eqref{cc13-1}, using \eqref{cc3} and \eqref{p14},  together with the Poincar\' e inequality, we deduce the final estimate:
\begin{align}
	&\intO{  \left[ \frac{\vt_s - \vt}{\vt}   \frac{\kappa (\vt) |\Grad \vt|^2 }{\vt} + \frac{\vt - \vt_s}{\vt_s} 
	\frac{ \kappa (\vt_s) |\Grad \vt_s|^2 }{\vt_s} + 
	\left(  \frac{\kappa (\vt) \Grad \vt }{\vt} - \frac{\kappa (\vt_s) \Grad \vt_s }{\vt_s}  \right) \cdot (\Grad \vt - \Grad \vt_s) \right] }  \br 
	&\gtrsim \| \vt - \vt_s \|^2_{W^{1,2}_0 (\Omega)}  + \intO{(\mathds{1}_{\Ov{\vt} \leq \vt  }   + \mathds{1}_{0< \vt \leq \underline\vt  }  )   \frac{\kappa(\vt)}{\vt^2}  |\Grad \vt|^2 } +  \intO{  \mathds{1}_{\vt \geq \underline{\vt}}  \frac{\kappa(\vt)}{\vt^2}  |\Grad \vt - \Grad \vt_s|^2 }.
	\label{cc14-0}
\end{align}

\end{enumerate}

\subsection{Estimates involving viscosity}\label{ev}

We recall the terms related to viscosity in \eqref{c5} and write them in the form:
\begin{align}\label{ev1}
& \intO{  \mathbb{S}(\vt,\Grad (\vu -\vu_s))  : \Grad(\vu - \vu_s) }	 +  \intO{ \frac{\vt_s -\vt}{\vt}\mathbb{S}(\vt,\Grad\vu): \Grad\vu   +\frac{\vt - \vt_s}{\vt_s} \mathbb{S}(\vt_s, \Grad \vu_s) : \Grad \vu_s  } \br
& = \intO{ \frac{\vt_{s}}{\vt} \mathbb{S}(\vt,\Grad (\vu -\vu_s))  : \Grad(\vu - \vu_s) }	\br
& \quad +  \intO{ \frac{\vt_s -\vt}{\vt}\Big( \mathbb{S}(\vt,\Grad\vu_{s}): \Grad\vu  + \mathbb{S}(\vt,\Grad\vu): \Grad\vu_{s} - 2 \mathbb{S}(\vt,\Grad\vu_{s}): \Grad\vu_{s}  \Big)} \br
&  \quad+  \intO{ \frac{\vt_s -\vt}{\vt} \mathbb{S}(\vt,\Grad\vu_{s}): \Grad\vu_{s}  - \frac{\vt_s -\vt}{\vt_{s}} \mathbb{S}(\vt_{s},\Grad\vu_{s}): \Grad\vu_{s} }.
\end{align}
Similarly to the preceding part, the above quantity decomposes into three regions: 

\begin{enumerate}

\item {\bf Estimates for small $\vt$.}
With the choice of $\underline \vt$ and $\Ov{\vt}$ as in \eqref{cc3}, it is rather straightforward to deduce that 
\begin{align}
	&\mathds{1}_{0 < \vt \leq \underline{\vt} } \frac{\vt_s -\vt}{\vt}\mathbb{S}(\vt,\Grad\vu): \Grad\vu   + 	\mathds{1}_{0 < \vt \leq \underline{\vt} }	 \frac{\vt - \vt_s}{\vt_s} \mathbb{S}(\vt_s, \Grad \vu_s) : \Grad \vu_s  \br 
	&\quad \geq \mathds{1}_{0 < \vt \leq \underline{\vt} }  \frac{\underline\vt}{2\vt} \mathbb{S}(\vt,\Grad\vu): \Grad\vu  - \mathds{1}_{0 < \vt \leq \underline{\vt} }	 \mathbb{S}(\vt_s, \Grad \vu_s) : \Grad \vu_s \br
	& \quad \geq \mathds{1}_{0 < \vt \leq \underline{\vt} }  \frac{\underline\vt}{2\vt} \mathbb{S}(\vt,\Grad\vu): \Grad\vu  - \mathds{1}_{0 < \vt \leq \underline{\vt} }	 c(\underline \vt)  |\Grad \vu_{s}|^{2}  (\vt - \vt_{s})^{2}.
\nonumber\end{align}
	Thus
\begin{align}\label{ev2}
	&\intO{ \mathds{1}_{0 < \vt \leq \underline{\vt} }   \mathbb{S}(\vt,\Grad (\vu -\vu_s))  : \Grad(\vu - \vu_s) }	\br
	& \quad  +  \intO{ \mathds{1}_{0 < \vt \leq \underline{\vt} }  \frac{\vt_s -\vt}{\vt}\mathbb{S}(\vt,\Grad\vu): \Grad\vu   +\mathds{1}_{0 < \vt \leq \underline{\vt} }  \frac{\vt - \vt_s}{\vt_s} \mathbb{S}(\vt_s, \Grad \vu_s) : \Grad \vu_s  } \br 
	&\geq \intO{ \mathds{1}_{0 < \vt \leq \underline{\vt} }   \mathbb{S}(\vt,\Grad (\vu -\vu_s))  : \Grad(\vu - \vu_s) } \br
	&\quad + \intO{ \mathds{1}_{0 < \vt \leq \underline{\vt} }   \frac{\underline\vt}{2\vt} \mathbb{S}(\vt,\Grad\vu): \Grad\vu }
	 - c(\underline \vt)  \|\Grad \vu_{s}\|_{L^{\infty}(\Omega)}^{2}\intO{ \mathds{1}_{0 < \vt \leq \underline{\vt}}  (\vt -\vt_{s})^{2}}.  
	\end{align}

\item {\bf Estimates for intermediate $\vt$.}
In this region, we use the right-hand side of \eqref{ev1} and deduce 
\begin{align}
&\mathds{1}_{\underline{\vt} \leq \vt \leq \Ov{\vt} } \frac{\vt_s -\vt}{\vt}\Big( \mathbb{S}(\vt,\Grad\vu_{s}): \Grad\vu  + \mathbb{S}(\vt,\Grad\vu): \Grad\vu_{s} - 2 \mathbb{S}(\vt,\Grad\vu_{s}): \Grad\vu_{s}  \Big)  \br
&\quad \leq \mathds{1}_{\underline{\vt} \leq \vt \leq \Ov{\vt} }  c(\underline \vt, \Ov{\vt}) |\vt_s -\vt| | \Grad \vu - \Grad\vu_{s}|  |\Grad \vu_{s} |  \br
&\quad \leq \mathds{1}_{\underline{\vt} \leq \vt \leq \Ov{\vt} }  c(\underline \vt, \Ov{\vt})  |\Grad \vu_{s} |  ( |\vt_s -\vt|^{2}  + | \Grad \vu - \Grad\vu_{s}|^{2}),
\nonumber \end{align}
 and 
\begin{align}
&\mathds{1}_{\underline{\vt} \leq \vt \leq \Ov{\vt} }\Big( \frac{\vt_s -\vt}{\vt} \mathbb{S}(\vt,\Grad\vu_{s}): \Grad\vu_{s}  - \frac{\vt_s -\vt}{\vt_{s}} \mathbb{S}(\vt_{s},\Grad\vu_{s}): \Grad\vu_{s}   \Big)   \leq \mathds{1}_{\underline{\vt} \leq \vt \leq \Ov{\vt} }  c(\underline \vt, \Ov{\vt}) |\vt_s -\vt|^{2}  |\Grad \vu_{s} |^{2} . 
\nonumber \end{align}
Under the condition  $\|\Grad \vu_{s}\|_{L^{\infty}(\Omega)} \lesssim \ep  \leq  \ep_{0}$ with $\ep_{0}$ sufficiently small, we finally deduce
\begin{align}\label{ev3}
 &\intO{\mathds{1}_{\underline{\vt} \leq \vt \leq \Ov{\vt} }  \frac{\vt_{s}}{\vt} \mathbb{S}(\vt,\Grad (\vu -\vu_s))  : \Grad(\vu - \vu_s) }	\br
& \quad +  \intO{\mathds{1}_{\underline{\vt} \leq \vt \leq \Ov{\vt} }  \frac{\vt_s -\vt}{\vt}\Big( \mathbb{S}(\vt,\Grad\vu_{s}): \Grad\vu  + \mathbb{S}(\vt,\Grad\vu): \Grad\vu_{s} - 2 \mathbb{S}(\vt,\Grad\vu_{s}): \Grad\vu_{s}  \Big)} \br
&  \quad+  \intO{\mathds{1}_{\underline{\vt} \leq \vt \leq \Ov{\vt} } \Big( \frac{\vt_s -\vt}{\vt} \mathbb{S}(\vt,\Grad\vu_{s}): \Grad\vu_{s}  - \frac{\vt_s -\vt}{\vt_{s}} \mathbb{S}(\vt_{s},\Grad\vu_{s}): \Grad\vu_{s}\Big) } \br 
& \geq \intO{\mathds{1}_{\underline{\vt} \leq \vt \leq \Ov{\vt} }  \frac{\vt_{s}}{2 \vt} \mathbb{S}(\vt,\Grad (\vu -\vu_s))  : \Grad(\vu - \vu_s) } -  c(\underline \vt, \Ov{\vt}) \|\Grad \vu_{s}\|_{L^{\infty}(\Omega)} \intO{\mathds{1}_{\underline{\vt} \leq \vt \leq \Ov{\vt} } |\vt_s -\vt|^{2}}  .
\end{align}

\item {\bf Estimates for large $\vt$.}
In this region, we also use the right-hand side of \eqref{ev1} and deduce 
\begin{align}
&\mathds{1}_{\vt \geq \Ov{\vt} } \frac{\vt_s -\vt}{\vt}\Big( \mathbb{S}(\vt,\Grad\vu_{s}): \Grad\vu  + \mathbb{S}(\vt,\Grad\vu): \Grad\vu_{s} - 2 \mathbb{S}(\vt,\Grad\vu_{s}): \Grad\vu_{s}  \Big)  \br
&\quad \leq \mathds{1}_{\vt \geq \Ov{\vt} } 2  \mu(\vt) | \Grad \vu - \Grad\vu_{s}|  |\Grad \vu_{s} |  \br 
&\quad \leq \mathds{1}_{\vt \geq \Ov{\vt} } |\Grad \vu_{s} |  2  \mu(\vt)^{2} + \mathds{1}_{\vt \geq \Ov{\vt} } |\Grad \vu_{s} |   | \Grad \vu - \Grad\vu_{s}|^{2}  .
\nonumber  \end{align}
 and 
\begin{align}
\mathds{1}_{\vt \geq \Ov{\vt} }  \Big( \frac{\vt_s -\vt}{\vt} \mathbb{S}(\vt,\Grad\vu_{s}): \Grad\vu_{s}  - \frac{\vt_s -\vt}{\vt_{s}} \mathbb{S}(\vt_{s},\Grad\vu_{s}): \Grad\vu_{s}   \Big)  \leq \mathds{1}_{\vt \geq \Ov{\vt} } 2  \mu(\vt) |  |\Grad \vu_{s} |^{2}.  
\nonumber \end{align}

In accordance with hypothesis \eqref{S12}, there holds
\begin{align}
 &\intO{\mathds{1}_{\vt \geq \Ov{\vt} } (\mu(\vt) + \mu(\vt)^{2} ) }	\leq c(\overline \vt)  \intO{\mathds{1}_{\vt \geq \Ov{\vt} } \kappa(\vt)}	  \leq c(\overline \vt)  \intO{ \mathds{1}_{\vt \geq \Ov{\vt} } \frac{\kappa (\vt)}{\vt^2} |\Grad \vt |^2 } .
  \nonumber \end{align}
Now, using again $\|\Grad \vu_{s}\|_{L^{\infty}(\Omega)} \lesssim \ep  \leq  \ep_{0}$ with $\ep_{0}$ sufficiently small, we finally deduce
\begin{align}\label{ev4}
 &\intO{\mathds{1}_{\vt \geq \Ov{\vt} }   \frac{\vt_{s}}{\vt} \mathbb{S}(\vt,\Grad (\vu -\vu_s))  : \Grad(\vu - \vu_s) }	\br
& \quad +  \intO{ \mathds{1}_{\vt \geq \Ov{\vt} }   \frac{\vt_s -\vt}{\vt}\Big( \mathbb{S}(\vt,\Grad\vu_{s}): \Grad\vu  + \mathbb{S}(\vt,\Grad\vu): \Grad\vu_{s} - 2 \mathbb{S}(\vt,\Grad\vu_{s}): \Grad\vu_{s}  \Big)} \br
&  \quad+  \intO{ \mathds{1}_{\vt \geq \Ov{\vt} } \Big( \frac{\vt_s -\vt}{\vt} \mathbb{S}(\vt,\Grad\vu_{s}): \Grad\vu_{s}  - \frac{\vt_s -\vt}{\vt_{s}} \mathbb{S}(\vt_{s},\Grad\vu_{s}): \Grad\vu_{s}\Big) } \br 
& \geq \intO{ \mathds{1}_{\vt \geq \Ov{\vt} }   \frac{\vt_{s}}{2 \vt} \mathbb{S}(\vt,\Grad (\vu -\vu_s))  : \Grad(\vu - \vu_s) } -  c( \Ov{\vt}) \|\Grad \vu_{s}\|_{L^{\infty}(\Omega)}  \intO{ \mathds{1}_{\vt \geq \Ov{\vt} } \frac{\kappa (\vt)}{\vt^2} |\Grad \vt |^2 } .
\end{align}

\item {\bf Conclusion.}

Putting the above estimates together, notably \eqref{ev2}, \eqref{ev3}, \eqref{ev4}, together with hypothesis \eqref{S12} on the viscosity  coefficients, Korn inequality and Poincar\'e inequality,  we deduce 
\begin{align}\label{ev5}
& \intO{  \mathbb{S}(\vt,\Grad (\vu -\vu_s))  : \Grad (\vu - \vu_s)  }+  \intO{ \frac{\vt_s -\vt}{\vt}\mathbb{S}(\vt,\Grad\vu): \Grad\vu   +\frac{\vt - \vt_s}{\vt_s} \mathbb{S}(\vt_s, \Grad \vu_s) : \Grad \vu_s  } \br
& \quad \gtrsim \| \vu - \vu_s \|^2_{W^{1,2}_0 (\Omega)}  + \intO{ \mathds{1}_{0 < \vt \leq \underline{\vt} }   \frac{1}{\vt} \mathbb{S}(\vt,\Grad\vu): \Grad\vu }  \br 
& \qquad - c(\underline \vt, \Ov{\vt}) \|\Grad \vu_{s}\|_{L^{\infty}(\Omega)} \Big(\|\vt_s -\vt\|_{L^{2}(\Omega)}^{2} + \intO{ \mathds{1}_{\vt \geq \Ov{\vt} } \frac{\kappa (\vt)}{\vt^2} |\Grad \vt |^2 } \Big)
\end{align}

Combining \eqref{ev5} with \eqref{cc14-0}, we get the final estimate concerning the terms involving dissipation: 
\begin{align}\label{cc14}
& \intO{  \mathbb{S}(\vt,\Grad (\vu -\vu_s))  : \Grad (\vu - \vu_s)  }+  \intO{ \frac{\vt_s -\vt}{\vt}\mathbb{S}(\vt,\Grad\vu): \Grad\vu   +\frac{\vt - \vt_s}{\vt_s} \mathbb{S}(\vt_s, \Grad \vu_s) : \Grad \vu_s  } \br 
& + \intO{  \left[ \frac{\vt_s - \vt}{\vt}   \frac{\kappa (\vt) |\Grad \vt|^2 }{\vt} + \frac{\vt - \vt_s}{\vt_s}  	\frac{ \kappa (\vt_s) |\Grad \vt_s|^2 }{\vt_s} +  	\left(  \frac{\kappa (\vt) \Grad \vt }{\vt} - \frac{\kappa (\vt_s) \Grad \vt_s }{\vt_s}  \right) \cdot (\Grad \vt - \Grad \vt_s) \right] }  \br 
&\gtrsim \| \vu - \vu_s \|^2_{W^{1,2}_0 (\Omega)}  +    \intO{ \mathds{1}_{0 < \vt \leq \underline{\vt} }   \frac{1}{\vt} \mathbb{S}(\vt,\Grad\vu): \Grad\vu }  + \| \vt - \vt_s \|^2_{W^{1,2}_0 (\Omega)}   \br 
& \quad + \intO{(\mathds{1}_{\Ov{\vt} \leq \vt  }   + \mathds{1}_{0< \vt \leq \underline\vt  }  )   \frac{\kappa(\vt)}{\vt^2}  |\Grad \vt|^2 } +  \intO{  \mathds{1}_{\vt \geq \underline{\vt}}  \frac{\kappa(\vt)}{\vt^2}  |\Grad \vt - \Grad \vt_s|^2 } . 
\end{align}

\end{enumerate}

\subsection{Pressure estimates and density damping}
\label{D}

A brief inspection of \eqref{cc14} reveals two ``damping'' terms, namely,
\[
 \| \vu - \vu_s \|^2_{W^{1,2}_0 (\Omega)} \ \mbox{and}\ \| \vt - \vt_s \|^2_{W^{1,2}_0 (\Omega)}
\]
acting on the velocity and temperature deviation, respectively. Obtaining similar damping on the density 
is quite delicate as ``mass dissipation'' is obviously not included in the NSF system. The estimates below seem therefore new and of independent 
interest.

We rewrite the pressure law \eqref{S9} in the form
\begin{equation} \label{c11}
	p(\vr, \vt) = p_\infty \vr^{\frac 53} + p_m(\vr, \vt) + \frac{a}{3} \vt^4,\quad  p_m(\vr, \vt) := \vt^{\frac{5}{2}} P_m \left( \frac{\vr}{\vt^{\frac{3}{2}}} \right).
\end{equation} 
Further, set $T_k(s) $ be a smooth non--decreasing function on $[0,\infty)$ such that 
\[ 
T_k (s) = s \ \mbox{if $0\leq s\leq k$}, \quad T_k (s) = k+1 \ \mbox{if $ s\geq k+2$}, \quad 0\leq T_{k}' \leq 1.
\]
Finally, we introduce the nowadays standard inverse of the divergence operator frequently called Bogovskii's operator $\mathcal{B}$, 
see Bogovskii \cite{BOG}, Galdi \cite{GALN}, Gei{\ss}ert,  Heck, and Hieber \cite{GEHEHI}: 
\begin{itemize}
	\item 
	$\mathcal{B}$ is a bounded linear operator mapping $L^q_0(\Omega)$ to $W^{1,q}_0(\Omega; R^3)$, where 
	$\Omega$ is a bounded Lipschitz domain in $R^3$ and the symbol $L^q_0$ denotes the space of $L^q-$integrable 
	functions with zero mean:	
	\begin{align} \label{AS6}
		\Div \mathcal{B}[f] = f \ \mbox{in}\ Q \ \mbox{for any}\ f \in L^q(\Omega), \ \langle f \rangle \equiv \frac{1}{|\Omega|} \intO{f} = 0, \
		\mathcal{B}[f]|_{\partial \Omega} = 0, \br 
		\| \mathcal{B}[f]  \|_{W^{1,q}(\Omega; R^3)} \lesssim \| f \|_{L^q(\Omega)},\ 1 < q < \infty.
	\end{align}
	\item If, moreover, $f = \Div \vc{g}$, $\vc{g} \in L^r(\Omega; R^3)$, $\vc{g} \cdot \vc{n}|_{\partial \Omega} = 0$, then 
	\begin{equation} \label{AS7}
		\| \mathcal{B}[\Div \vc{g} ] \|_{L^r(\Omega; R^3)} \lesssim \| \vc{g} \|_{L^r(\Omega; R^3)},\ 1 < r < \infty.
	\end{equation}

\end{itemize}
 
The desired density estimates are obtained by testing the momentum equation \eqref{w3} on
\[
\phiB:=\mathcal{B} \left[ T_k(\vr) - T_k(\vr_s) - \langle T_k(\vr) - T_k(\vr_s) \rangle \right].
\]
To this end, we first rewrite \eqref{w3} in the form,
\begin{align}
	- \frac{\D}{\dt} \intO{ \vr \vu \cdot \phiB } &+ \intO{ \Big( \vr \vu \cdot \partial_t \phiB + \vr \vu \otimes \vu : \Grad \phiB + 
		p(\vr, \vt) \Div \phiB \Big) } \br &= \intO{ \mathbb{S}(\vt, \Grad \vu) : \Grad \phiB } -
	\intO{ \vr \Grad G \cdot \phiB }   	
	\label{w3A}
\end{align}
in $\mathcal{D}'(T, \infty)$.

Now, a direct computation yields
\begin{align}
	 &\intO{ \Big( p(\vr, \vt_s) - p(\vr_s, \vt_s) \Big) \Big(T_{k}(\vr) - T_{k}(\vr_s)\Big) }	\br \quad &= 
\intO{ \Big( p(\vr, \vt) - p(\vr_s, \vt_s) \Big) \Big(T_{k}(\vr) - T_{k}(\vr_s)\Big)} + \intO{ \Big( p(\vr, \vt_s) - p(\vr, \vt) \Big) \Big(T_{k}(\vr) - T_{k}(\vr_s)\Big)}	\br \quad &=  \intO{ \Big( p(\vr, \vt) - p(\vr_s, \vt_s) \Big) \Div \phiB } + \langle T_k(\vr) - T_k(\vr_s) \rangle \intO{ \Big( p(\vr, \vt) - p(\vr_s, \vt_s) \Big)  }    \br
 &\quad + \intO{ \Big( p_m(\vr, \vt_s) - p_m(\vr, \vt) \Big) \Big(T_{k}(\vr) - T_{k}(\vr_s)\Big) }  + \frac{a}{3} \intO{ \Big( \vt_s^4 - \vt^4 \Big) \Big(T_{k}(\vr) - T_{k}(\vr_s)\Big) } \br
 & = \intO{ \Big( p(\vr, \vt) - p(\vr_s, \vt_s) \Big) \Div \phiB } + \sum_{i=1}^{3}I_{i}.
	\label{c12}
	\end{align}
Furthermore, in accordance with the weak formulation \eqref{w3A} and the stationary equation \eqref{p10}, 
	\begin{align}
	&  \intO{ \Big( p(\vr, \vt) - p(\vr_s, \vt_s) \Big) \Div \phiB } \br
&= \intO{ \Big(\mathbb{S} (\vt, \Grad \vu)  - \mathbb{S} (\vt_{s}, \Grad \vu_{s})  \Big): \Grad  \phiB} - 
\intO{\Big( \vr \vu \otimes \vu - \vr_{s} \vu_{s} \otimes \vu_{s} \Big) : \Grad \phiB}  \br 
& \quad - \intO{ (\vr-\vr_{s})\Grad G \cdot \phiB }  + \intO{ \Big(\vr \vu - \vr_{s} \vu_{s} \Big) \cdot \partial_{t} \phiB} - \frac{\D}{\dt} \intO{ \Big(\vr \vu -\vr_{s} \vu_{s} \Big)\cdot \phiB} \br 
 &= \sum_{i=4}^8 I_i.
	\label{c12-1}
	\end{align}

\medskip

We write
\begin{align} 
&\intO{ \Big( p(\vr, \vt_s) - p(\vr_s, \vt_s) \Big) \Big(T_{k}(\vr) - T_{k}(\vr_s)\Big) }\br
& = \intO{ \mathds{1}_{ \underline \vr \leq \vr \leq \overline \vr} \Big( p(\vr, \vt_s) - p(\vr_s, \vt_s) \Big) \Big(T_{k}(\vr) - T_{k}(\vr_s)\Big) }\br
 & \quad + \intO{ \mathds{1}_{ \vr \geq \overline \vr} \Big( p(\vr, \vt_s) - p(\vr_s, \vt_s) \Big) \Big(T_{k}(\vr) - T_{k}(\vr_s)\Big) } \br 
 & \quad + \intO{ \mathds{1}_{ \vr \leq \underline \vr } \Big( p(\vr, \vt_s) - p(\vr_s, \vt_s) \Big) \Big(T_{k}(\vr) - T_{k}(\vr_s)\Big) } ,
 \label{c13-1}
\end{align} 
where we have set positive constants $k$, $\underline \vr$ and $\overline \vr$ so that
\begin{equation}\label{def-bar-vr}
0< \underline \vr \leq \frac{1}{2}\inf_{\Omega}  \vr_{s} < 2 \sup_{\Omega}  \vr_{s} \leq \overline \vr \leq k. 
\end{equation}
 As the pressure $P\in C^{1}[0, \infty)$ is a strictly increasing function with $P'(Z) > 0$ for all $Z \geq 0$ (see \eqref{S5}), there holds for any $0< Z_{1} < Z_{2} < \infty$ that
$$
P'(Z) \gtrsim 1, \quad \mbox{for all} \ Z \in [Z_{1}, Z_{2}].
$$
Consequently, 
\begin{align} 
& \intO{ \mathds{1}_{ \underline \vr \leq \vr \leq \overline \vr} \Big( p(\vr, \vt_s) - p(\vr_s, \vt_s) \Big) \Big(T_{k}(\vr) - T_{k}(\vr_s)\Big) }\br
 & = \intO{ \mathds{1}_{ \underline \vr \leq \vr \leq \overline \vr}  \vt_{s} P'\big(\frac{\tilde \vr}{\vt_{s}^{\frac 32}}\big)(\vr - \vr_{s}) \Big(T_{k}(\vr) - T_{k}(\vr_s)\Big) }\br
 & \gtrsim \intO{ \mathds{1}_{ \underline \vr \leq \vr \leq \overline \vr}   (\vr - \vr_{s})^{2}  },
 \label{c13-2}
\end{align} 
where $\tilde \vr \in [\underline \vr , \overline \vr]$.

Using the specific form of $p$ in \eqref{c11} and the sublinear growth of $P_{m}$ in \eqref{pm-infty-2}, and choosing $\overline \vr$ suitably large, implies
\begin{align} 
& \intO{ \mathds{1}_{ \vr \geq \overline \vr} \Big( p(\vr, \vt_s) - p(\vr_s, \vt_s) \Big) \Big(T_{k}(\vr) - T_{k}(\vr_s)\Big) }\br
 & = \intO{  \mathds{1}_{ \vr \geq \overline \vr}  \Big[p_{\infty} (\vr^{\frac 53} - \vr_{s}^{\frac 53}) +   \vt_{s}^{\frac 52} \Big( P_{m}\big(\frac{ \vr}{\vt_{s}^{\frac 32}}\big) -P_{m}\big(\frac{\vr_{s}}{\vt_{s}^{\frac 32}}\big) \Big)\Big] \Big(T_{k}(\vr) - T_{k}(\vr_s)\Big) }\br
 & \gtrsim \intO{ \mathds{1}_{ \vr \geq \overline \vr}   \vr ^{\frac 53}  }.
 \label{c13-3}
\end{align} 

As $P_{m}(0) = P(0) = 0$ and $P_{m}(Z)>0$ for all $Z>0$, by choosing $\underline \vr$ suitably small there holds 
$$
P_{m}\big(\frac{ \vr}{\vt_{s}^{\frac 32}}\big) < P_{m}\big(\frac{\vr_{s}}{\vt_{s}^{\frac 32}}\big), \quad \mbox{for all $0\leq \vr \leq \underline \vr$}.
$$ 
Thus,
\begin{align} 
& \intO{ \mathds{1}_{ \vr \leq \underline \vr} \Big( p(\vr, \vt_s) - p(\vr_s, \vt_s) \Big) \Big(T_{k}(\vr) - T_{k}(\vr_s)\Big) }\br
 & = p_{\infty}\intO{  \mathds{1}_{ \vr \leq \underline \vr}   (\vr^{\frac 53} - \vr_{s}^{\frac 53}) (\vr  - \vr_s) } +\intO{  \mathds{1}_{ \vr \leq \underline \vr}   \vt_{s}^{\frac 52} \Big( P_{m}\big(\frac{ \vr}{\vt_{s}^{\frac 32}}\big) -P_{m}\big(\frac{\vr_{s}}{\vt_{s}^{\frac 32}}\big) \Big) (\vr  - \vr_s) }\br
 & \geq p_{\infty}  \intO{  \mathds{1}_{ \vr \leq \underline \vr}  (\vr^{\frac 53} - \vr_{s}^{\frac 53}) (\vr  - \vr_s) } \br
 &  \gtrsim\intO{  \mathds{1}_{ \vr \leq \underline \vr} }.
 \label{c13-4}
\end{align} 

Summing up the estimates in \eqref{c13-1}--\eqref{c13-4} finally gives
\begin{align} 
& \intO{  \Big( p(\vr, \vt_s) - p(\vr_s, \vt_s) \Big) \Big(T_{k}(\vr) - T_{k}(\vr_s)\Big) }\br
 &  \gtrsim \intO{ \mathds{1}_{ \underline \vr \leq \vr \leq \overline \vr}   (\vr - \vr_{s})^{2}  } + \intO{ \mathds{1}_{ \vr \geq \overline \vr}   \vr ^{\frac 53}  } +  \intO{  \mathds{1}_{ \vr \leq \underline \vr} } .
 \label{c13}
\end{align}

\medskip

Next, we estimate the integrals $I_1-I_8$ in terms of dissipation. We start with $I_1, I_2, I_3$ in \eqref{c12}.
\begin{itemize}
	\item {\bf Integral $I_1$.} First,
\begin{align} 
I_{1} &  = \frac{1}{|\Omega|} \intO{ T_k(\vr) - T_k(\vr_s) }  \intO{ \Big( p(\vr, \vt) - p(\vr_s, \vt_s) \Big)  }  \br
& = \frac{1}{|\Omega|} \intO{ T_k(\vr) - T_k(\vr_s) }  \intO{ \Big( p(\vr, \vt) - p(\vr_s, \vt) \Big)  + \Big( p(\vr_{s}, \vt) - p(\vr_s, \vt_s) \Big) }  .
\label{I1}
\end{align}
Seeing that  $T_k(\vr_s) = \vr_s$ as soon as $k$ satisfies \eqref{def-bar-vr} we obtain
\begin{equation} \label{I1-0}
\intO{ \Big( T_k(\vr) - T_k(\vr_s) \Big) } = \intO{ \Big(  T_k(\vr) - \vr + \vr - \vr_s \Big) } = 
\intO{ \Big( T_k(\vr) - \vr \Big) } \leq 0.
\end{equation}

In accordance with hypothesis \eqref{convex-pm}, the function $P$ is convex; whence 
$\vr \mapsto p(\vr, \vt)$ is convex for any fixed $\vt$. In particular, 
\[
p(\vr, \vt) - p(\vr_s, \vt) \geq \frac{\partial p}{\partial \vr}(\vr_{s}, \vt) (\vr - \vr_s).
\]
Consequently, in view of \eqref{I1-0},
\begin{align} 
	& \intO{ \Big( T_k(\vr) - T_k(\vr_s) \Big) }   \intO{ \Big( p(\vr, \vt) - p(\vr_s, \vt) \Big) } \br 
	&\quad \leq \intO{  \frac{\partial p}{\partial \vr}(\vr_{s}, \vt)(\vr - \vr_s) } \intO{ \Big( T_k(\vr) - T_k(\vr_s) \Big) } \br 
	&\quad  = \intO{  \frac{\partial p}{\partial \vr}(\vr_{s}, \vt_{s}) (\vr - \vr_s)  + \Big(\frac{\partial p}{\partial \vr}(\vr_{s}, \vt) - \frac{\partial p}{\partial \vr}(\vr_{s}, \vt_{s}) \Big)(\vr - \vr_s)  } \intO{ \Big( T_k(\vr) - T_k(\vr_s) \Big) }. 
 	\nonumber
\end{align}	
Moreover, as
\[
\intO{ (\vr - \vr_s) } = 0,
\]
we get
\begin{align} 
	& \intO{ \Big( T_k(\vr) - T_k(\vr_s) \Big) }   \intO{ \Big( p(\vr, \vt) - p(\vr_s, \vt) \Big) } \br
		& \leq  
		\intO{  \Big( \frac{\partial p}{\partial \vr}(\vr_{s}, \vt_{s})  -  \frac{\partial p}{\partial \vr}(\vr_{c}, \vt_{c})  \Big) (\vr - \vr_s)  } \intO{ \Big( T_k(\vr) - T_k(\vr_s) \Big) } \br
		&+ \intO{  \Big( \frac{\partial p}{\partial \vr}(\vr_{s}, \vt) - \frac{\partial p}{\partial \vr}(\vr_{s}, \vt_{s}) \Big)(\vr - \vr_s)  } \intO{ \Big( T_k(\vr) - T_k(\vr_s) \Big) },
	\label{I1-1} 	
\end{align}
for arbitrary constants $\vr_{c}  > 0, \ \vt_{c}>0$. However, as $\vr_s$ and $\vt_{s}$ are close to constants in the sense of \eqref{p14},  we get 
\[
\Big| \frac{\partial p}{\partial \vr}(\vr_{s}, \vt_{s})  -  \frac{\partial p}{\partial \vr}(\vr_{c}, \vt_{c})  \Big | \lesssim \ep 
\]
for suitable constants $\vr_{c}$ and $\vt_{c}$. Thus,
\begin{align} 
	& \intO{  \Big( \frac{\partial p}{\partial \vr}(\vr_{s}, \vt_{s})  -  \frac{\partial p}{\partial \vr}(\vr_{c}, \vt_{c})  \Big) (\vr - \vr_s)  } \intO{ \Big( T_k(\vr) - T_k(\vr_s) \Big) } \br
	& \quad \lesssim \ep  \intO{  |\vr - \vr_{s}| } \intO{ \vr - T_{k} (\vr) } \lesssim \ep  \intO{ \mathds{1}_{\vr \geq \overline \vr} \vr  } \lesssim  \ep  \intO{ \mathds{1}_{\vr \geq \overline \vr} \vr^{\frac{5}{3}}  },
	\nonumber
\end{align}
which can be absorbed by the integrals on the right--hand side of \eqref{c13} provided $\ep$ is sufficiently small. 

\medskip

For the other integral in \eqref{I1-1}, by \eqref{I1-0} and the fact $P_{m}\in C^{1}[0, \infty)$ with $P_{m}'$ satisfying  \eqref{pm-infty-3}, we have
\begin{align} 
	& \intO{  \Big( \frac{\partial p}{\partial \vr}(\vr_{s}, \vt) - \frac{\partial p}{\partial \vr}(\vr_{s}, \vt_{s}) \Big)(\vr - \vr_s)  } \intO{ \Big( T_k(\vr) - T_k(\vr_s) \Big) } \br
	& \quad =  \intO{  \Big( \vt P_{m}'\big(\frac{\vr_{s}}{\vt^{\frac{3}{2}}}\big)  - \vt_{s} P_{m}' \big(\frac{\vr_{s}}{\vt_{s}^{\frac{3}{2}}}\big)  \Big)(\vr - \vr_s)  }  \intO{ \Big( T_k(\vr) - T_k(\vr_s) \Big) } \br
	& \quad \lesssim   \intO{ \mathds{1}_{ \underline\vt \leq \vt \leq \overline  \vt} | \vt - \vt_{s} | | \vr - \vr_s |  }  \intO{ \mathds{1}_{\vr \geq \overline \vr} \vr }  +  \intO{ \mathds{1}_{ \vt \geq \overline  \vt} \vt | \vr - \vr_s |  }   \intO{ \mathds{1}_{\vr \geq \overline \vr} \vr  } \br 
	& \qquad  +  \intO{ \mathds{1}_{ \vt \leq \underline  \vt}   | \vr - \vr_s |  }   \intO{ \mathds{1}_{\vr \geq \overline \vr} \vr  }.
		\label{I1-3} 	
\end{align}
For the integrals on the right-hand side of \eqref{I1-3}, 
\begin{align} 
	 & \intO{ \mathds{1}_{ \underline\vt \leq \vt \leq \overline  \vt} | \vt - \vt_{s} | | \vr - \vr_s |  }  \intO{ \mathds{1}_{\vr \geq \overline \vr} \vr } \br 
	 & \quad   \leq \|\vt-\vt_{s}\|_{L^{6}(\Omega)} \| \vr -\vr_{s}\|_{L^{\frac 65}(\Omega)}   \intO{ \mathds{1}_{\vr \geq \overline \vr} \vr }  \br 
	  & \quad \leq C(\delta)   \|\vt-\vt_{s}\|_{W^{1,2}(\Omega)}^{2} + \delta\intO{ \mathds{1}_{\vr \geq \overline \vr} \vr^{\frac 53}  } ,
	\label{I1-4} 	
\end{align}
\begin{align} 
		&  \intO{ \mathds{1}_{ \vt \geq \overline  \vt} \vt | \vr - \vr_s |  }   \intO{ \mathds{1}_{\vr \geq \overline \vr} \vr  } \br 
		&\quad  \leq  C(\delta)  \intO{ \mathds{1}_{ \vt \geq \overline  \vt} \vt^{\frac{5}{2}} }   \left(\intO{  | \vr - \vr_s |^{\frac 53}  } \right)^{\frac 32}  +  \delta \left(\intO{ \mathds{1}_{\vr \geq \overline \vr} \vr  } \right)\br
		& \quad \leq C(\delta) \intO{ \mathds{1}_{\vt \geq \Ov{\vt} } \frac{\kappa (\vt)}{\vt^2} |\Grad \vt |^2 } + \delta \intO{ \mathds{1}_{\vr \geq \overline \vr} \vr^{\frac 53}  } ,
		\label{I1-5} 	
\end{align}
 and
\begin{align} 
	& \intO{ \mathds{1}_{ \vt \leq \underline  \vt}  | \vr - \vr_s |  }   \intO{ \mathds{1}_{\vr \geq \overline \vr} \vr  } \br 
	& \quad \leq  C(\underline  \vt, \overline \vt)  \intO{ \mathds{1}_{ \vt \leq  \underline \vt}  |\vt-\vt_{s}|  | \vr - \vr_s |  }  \intO{ \mathds{1}_{\vr \geq \overline \vr} \vr  } \br
	& \quad   \leq  C(\underline  \vt, \overline \vt)  \|\vt-\vt_{s}\|_{L^{6}(\Omega)} \| \vr -\vr_{s}\|_{L^{\frac 65}(\Omega)}   \intO{ \mathds{1}_{\vr \geq \overline \vr} \vr }  \br 
	  & \quad \leq C(\delta)   \|\vt-\vt_{s}\|_{W^{1,2}(\Omega)}^{2} + \delta\intO{ \mathds{1}_{\vr \geq \overline \vr} \vr^{\frac 53}  } ,
\label{I1-6} 	
\end{align}
where we used the estimates in \eqref{cc11}--\eqref{cc12}, together with the essential fact that the norm $\| \vr \|_{L^{\frac{5}{3}}(\Omega)}$ is uniformly bounded for large times thanks to the existence of the  bounded absorbing set, cf. \eqref{A6}. Consequently, the integral is controlled by dissipation. 

\medskip

Seeing \eqref{I1}, to finish the estimates of $I_{1}$ , we are left with the integral
\begin{align} 
 \intO{ T_k(\vr) - T_k(\vr_s) }  \intO{ \Big( p(\vr_{s}, \vt) - p(\vr_s, \vt_s) \Big) }  .
\nonumber 
\end{align}
 By the properties of $P$ and $P_{m}$ in \eqref{S5}, \eqref{pm-infty-2} and \eqref{pm-infty-3}, 
\begin{align} 
 \intO{ \big| p(\vr_{s}, \vt) - p(\vr_s, \vt_s) \big |} \leq   C(\underline \vt, \overline \vt) \left( \intO{  \mathds{1}_{ \underline \vt \leq \vt \leq \overline \vt} | \vt -  \vt_s | } +  \intO{  \mathds{1}_{ \vt \geq \overline \vt  } \vt^{4}  }  + \intO{  \mathds{1}_{ \vt \leq \underline \vt  }  } \right).
\nonumber 
\end{align}
Similarly as in \eqref{I1-4}--\eqref{I1-6},
  \begin{align} 
	  \intO{ \mathds{1}_{ \underline\vt \leq \vt \leq \overline  \vt} | \vt - \vt_{s} |  }  \intO{ \mathds{1}_{\vr \geq \overline \vr} \vr }  \leq C(\delta)   \|\vt-\vt_{s}\|_{W^{1,2}(\Omega)}^{2} + \delta\intO{ \mathds{1}_{\vr \geq \overline \vr} \vr^{\frac 53}  } ,
\nonumber	
\end{align}
\begin{align} 
		  \intO{ \mathds{1}_{ \vt \geq \overline  \vt} \vt^{4}  }   \intO{ \mathds{1}_{\vr \geq \overline \vr} \vr  }   \lesssim  \intO{ \mathds{1}_{ \vt \geq \overline  \vt} \vt^{4} } \lesssim \intO{ \mathds{1}_{\vt \geq \Ov{\vt} } \frac{\kappa (\vt)}{\vt^2} |\Grad \vt |^2 } ,
	\nonumber 
\end{align}
and
\begin{align} 
	 \intO{ \mathds{1}_{ \vt \leq \underline  \vt}   }   \intO{ \mathds{1}_{\vr \geq \overline \vr} \vr  } 
	&  \leq  C(\underline  \vt, \overline \vt)  \intO{ \mathds{1}_{ \vt \leq  \underline \vt}  |\vt-\vt_{s}|  }  \intO{ \mathds{1}_{\vr \geq \overline \vr} \vr  } \br
	  &  \leq C(\delta)   \|\vt-\vt_{s}\|_{W^{1,2}(\Omega)}^{2} + \delta\intO{ \mathds{1}_{\vr \geq \overline \vr} \vr^{\frac 53}  } .
\nonumber
\end{align}

Summing up the above estimates implies 
\begin{align} 
	I_{1} \leq C(\delta)  \|\vt-\vt_{s}\|_{W^{1,2}(\Omega)}^{2} +  C(\delta)\intO{ \mathds{1}_{\vt \geq \Ov{\vt} } \frac{\kappa (\vt)}{\vt^2} |\Grad \vt |^2 }  + \delta\intO{ \mathds{1}_{\vr \geq \overline \vr} \vr^{\frac 53}  }.
\label{I1-f} 	
\end{align}
Here and hereafter,  $\delta$ denotes an arbitrary positive constant which will be chosen suitably small and $C(\delta)$ is a constant depending on $\delta.$

\item {\bf Integral $I_2$.} Next we consider 
\begin{align} 
	I_{2} =  \intO{ \Big( p_m(\vr, \vt_s) - p_m(\vr, \vt) \Big) \Big(T_{k}(\vr) - T_{k}(\vr_s)\Big) } .
\nonumber 
\end{align}

We shall decompose integral in $I_2$ into different regions with respect to the range of $\vt$ and $\vr$. We start with the region where $\vt \geq \overline \vt$. Note that $p_m(\vr, \vt)$ is increasing in $\vt \in (0,\infty)$ due to the fact
\begin{align} 
\frac{\partial p_m}{\partial \vt} (\vr, \vt)  = 3\vt^{\frac 32} \Big(\frac{5}{3} P_m(Z) - P_m'(Z) Z \Big) \geq 0, \ Z = \frac{\vr}{\vt^\frac 32}.
\label{I2-2} 	
\end{align}
Thus
\begin{align} 
	I_{2,1} =  \intO{ \mathds{1}_{\vt\geq \overline \vt }  \mathds{1}_{ \vr \geq  \overline \vr}\Big( p_m(\vr, \vt_s) - p_m(\vr, \vt) \Big) \Big(T_{k}(\vr) - T_{k}(\vr_s)\Big) } \leq 0.
\nonumber 
\end{align}

The fact $P_m\in C^{1}[0, \infty)$ implies
\begin{align} 
	I_{2,2} & =  \intO{ \mathds{1}_{\vt\geq \overline \vt }  \mathds{1}_{ \vr \leq  \overline \vr}\Big( p_m(\vr, \vt_s) - p_m(\vr, \vt) \Big) \Big(T_{k}(\vr) - T_{k}(\vr_s)\Big) } \br 
	& \lesssim \intO{ \mathds{1}_{\vt\geq \overline \vt }  \mathds{1}_{ \vr \leq  \overline \vr} \vt^\frac{5}{2}} \lesssim \intO{ \mathds{1}_{\vt\geq \overline \vt }   \frac{\kappa (\vt)}{\vt^2} |\Grad \vt |^2 }.
\nonumber 
\end{align}

We turn to the region $\vt \leq \underline \vt$. By virtue of \eqref{pm-infty-2}, 
\begin{align} 
	I_{2,3} & =  \intO{ \mathds{1}_{\vt\leq \underline \vt }  \mathds{1}_{ \vr \geq  \overline \vr}\Big( p_m(\vr, \vt_s) - p_m(\vr, \vt) \Big) \Big(T_{k}(\vr) - T_{k}(\vr_s)\Big) } \br 
	& \lesssim  \intO{ \mathds{1}_{\vt\leq \underline \vt }  \mathds{1}_{ \vr \geq  \overline \vr}\Big( \vt_{s}^{\frac 52} \big(\frac{\vr}{\vt_{s}^{\frac 32}}\big) + \vt^{\frac 52} \big(\frac{\vr}{\vt^{\frac 32}}\big)\Big) (k+1) } \br 
	& \lesssim  \intO{ \mathds{1}_{\vt\leq \underline \vt }  \mathds{1}_{ \vr \geq  \overline \vr} \vr | \vt - \vt_s|^{\frac 45} } \br
	& \lesssim  \intO{ \mathds{1}_{\vt\leq \underline \vt }  \mathds{1}_{ \vr \geq  \overline \vr}  \big( \delta \vr^\frac{5}{3} + C(\delta)| \vt - \vt_s|^2 \big)}\br 
	& \leq  C(\delta) \intO{  | \vt - \vt_s|^2  } +  \delta\intO{   \mathds{1}_{ \vr \geq  \overline \vr}   \vr^\frac{5}{3} }.
\nonumber 
\end{align}

Next, again by the fact $P_{m} \in C^{1}[0, \infty)$ and \eqref{pm-infty-2}, there holds
\begin{align} 
	I_{2,4} & =  \intO{ \mathds{1}_{\vt\leq \underline \vt }  \mathds{1}_{ \vr \leq  \overline \vr}\Big( p_m(\vr, \vt_s) - p_m(\vr, \vt) \Big) \Big(T_{k}(\vr) - T_{k}(\vr_s)\Big) } \br 
	& \lesssim  \intO{ \mathds{1}_{\vt\leq \underline \vt }  \mathds{1}_{ \vr \leq  \overline \vr}  \Big( \vt_{s}^{\frac 52} + \vt_{s}^{\frac 52} \big(\frac{\vr}{\vt_{s}^{\frac 32}}\big) + \vt^{\frac 52} + \vt^{\frac 52} \big(\frac{\vr}{\vt^{\frac 32}}\big)\Big) (k+1)} \br 
	& \lesssim  \intO{ \mathds{1}_{\vt\leq \underline \vt } | \vt - \vt_s|^2  }.
\nonumber 
\end{align}

Finally, we consider the region $\underline \vt \leq \vt \leq \overline \vt$: 
\begin{align} 
	I_{2,5} & =  \intO{ \mathds{1}_{\underline \vt \leq \vt \leq \overline \vt }  \mathds{1}_{ \underline \vr \leq \vr \leq \overline \vr}\Big( p_m(\vr, \vt_s) - p_m(\vr, \vt) \Big) \Big(T_{k}(\vr) - T_{k}(\vr_s)\Big) } \br 
	& \lesssim  \intO{ \mathds{1}_{\underline \vt \leq \vt \leq \overline \vt }  \mathds{1}_{ \underline \vr \leq \vr \leq \overline \vr} |\vt-\vt_s| |\vr-\vr_s| } \br 
	& \leq    C(\delta)\intO{ \mathds{1}_{\underline \vt \leq \vt \leq \overline \vt }   |\vt-\vt_s|^2 } + \delta \intO{   \mathds{1}_{ \underline \vr \leq \vr \leq \overline \vr}   |\vr-\vr_s|^2 }.
\nonumber 
\end{align}
Again by \eqref{pm-infty-2} and \eqref{pm-infty-3},  applying mean value theorem and using \eqref{I2-2} gives
\begin{align} 
	I_{2,6} & =  \intO{ \mathds{1}_{\underline \vt \leq \vt \leq \overline \vt }  \mathds{1}_{  \vr \leq \underline \vr}\Big( p_m(\vr, \vt_s) - p_m(\vr, \vt) \Big) \Big(T_{k}(\vr) - T_{k}(\vr_s)\Big) } \br 
	& \leq  \intO{ \mathds{1}_{\underline \vt \leq \vt \leq \overline \vt }  \mathds{1}_{  \vr \leq \underline \vr} \frac{\partial p_m}{\partial \vt} (\vr, \tilde \vt)  |\vt-\vt_s| (k+1)} \br 
	& \leq    C(\delta)\intO{ \mathds{1}_{\underline \vt \leq \vt \leq \overline \vt }   |\vt-\vt_s|^2 } + \delta \intO{   \mathds{1}_{ \vr \leq \underline \vr} },
\nonumber 
\end{align}
and
\begin{align} 
	I_{2,7} & =  \intO{ \mathds{1}_{\underline \vt \leq \vt \leq \overline \vt }  \mathds{1}_{  \vr \geq \overline \vr}\Big( p_m(\vr, \vt_s) - p_m(\vr, \vt) \Big) \Big(T_{k}(\vr) - T_{k}(\vr_s)\Big) } \br 
	& \leq  \intO{ \mathds{1}_{\underline \vt \leq \vt \leq \overline \vt }  \mathds{1}_{   \vr \geq \overline \vr} \frac{\partial p_m}{\partial \vt} (\vr, \hat \vt)  |\vt-\vt_s|  (k+1) } \br 
	& \lesssim  \intO{ \mathds{1}_{\underline \vt \leq \vt \leq \overline \vt }  \mathds{1}_{   \vr \geq \overline \vr} \vr |\vt-\vt_s|^{\frac 45} } \br 
	& \lesssim     \intO{ \mathds{1}_{\underline \vt \leq \vt \leq \overline \vt }  \mathds{1}_{   \vr \geq \overline \vr}  \big( \delta \vr^\frac{5}{3} + C(\delta)| \vt - \vt_s|^2 \big) } \br
	   & \leq  C(\delta) \intO{  | \vt - \vt_s|^2  } +  \delta\intO{   \mathds{1}_{ \vr \geq  \overline \vr}   \vr^\frac{5}{3} },
\nonumber 
\end{align}
with $\tilde \vt , \hat \vt \in [\underline \vt, \overline \vt]$. We have therefore completed the estimate of $I_2$:
\begin{align} 
	I_{2} & \leq  C(\delta )\intO{ \mathds{1}_{\vt\geq \overline \vt }   \frac{\kappa (\vt)}{\vt^2} |\Grad \vt |^2 } +  C(\delta)\|\vt-\vt_s\|_{W^{1,2}(\Omega)}^2 \br
	&\quad + \delta \intO{   \mathds{1}_{ \underline \vr \leq \vr \leq \overline \vr}   |\vr-\vr_s|^2 } + \delta \intO{   \mathds{1}_{ \vr \leq \underline \vr} }  +\delta \intO{   \mathds{1}_{ \vr \geq \overline \vr} \vr^\frac 53 }. 
\label{I2-f} 	
\end{align}

\item {\bf Integral $I_3$.} Now, we consider 
\begin{align} 
	I_{3} =  \frac{a}{3} \intO{ \Big( \vt_s^4 - \vt^4 \Big) \Big(T_{k}(\vr) - T_{k}(\vr_s)\Big) } .
\nonumber 
\end{align}
By arguments similar to the estimate of $I_2$, we have 
\begin{align} 
	I_{3,1} & =  \frac{a}{3} \intO{ \mathds{1}_{\vt \geq \overline \vt} \Big( \vt_s^4 - \vt^4 \Big) \Big(T_{k}(\vr) - T_{k}(\vr_s)\Big) } \br 
	& \lesssim  \intO{ \mathds{1}_{\vt \geq \overline \vt}  \vt^4 } \lesssim \intO{ \mathds{1}_{\vt\geq \overline \vt }   \frac{\kappa (\vt)}{\vt^2} |\Grad \vt |^2 } ,
\nonumber 
\end{align}
\begin{align} 
	I_{3,2} & =  \frac{a}{3} \intO{ \mathds{1}_{\vt \leq \underline \vt} \Big( \vt_s^4 - \vt^4 \Big) \Big(T_{k}(\vr) - T_{k}(\vr_s)\Big) } \lesssim  \intO{ \mathds{1}_{\vt \leq \underline \vt} |\vt - \vt_s|^2  } ,
\nonumber 
\end{align}
\begin{align} 
	I_{3,3} & =  \frac{a}{3} \intO{ \mathds{1}_{\underline \vt \leq \vt \leq \overline \vt }  \mathds{1}_{ \underline \vr \leq \vr \leq \overline \vr} \Big( \vt_s^4 - \vt^4 \Big) \Big(T_{k}(\vr) - T_{k}(\vr_s)\Big) } \br 
	& \lesssim  \intO{ \mathds{1}_{\underline \vt \leq \vt \leq \overline \vt }  \mathds{1}_{ \underline \vr \leq \vr \leq \overline \vr} |\vt-\vt_s| |\vr-\vr_s| } \br 
	& \leq    C(\delta)\intO{ \mathds{1}_{\underline \vt \leq \vt \leq \overline \vt }   |\vt-\vt_s|^2 } + \delta \intO{   \mathds{1}_{ \underline \vr \leq \vr \leq \overline \vr}   |\vr-\vr_s|^2 },
\nonumber 
\end{align}
\begin{align} 
	I_{3,4} & =  \intO{ \mathds{1}_{\underline \vt \leq \vt \leq \overline \vt }  \mathds{1}_{  \vr \leq \underline \vr}\Big( \vt_s^4 - \vt^4 \Big) \Big(T_{k}(\vr) - T_{k}(\vr_s)\Big) } \br 
	& \lesssim  \intO{ \mathds{1}_{\underline \vt \leq \vt \leq \overline \vt }  \mathds{1}_{  \vr \leq \underline \vr} |\vt-\vt_s|} \br 
	& \leq    C(\delta)\intO{ \mathds{1}_{\underline \vt \leq \vt \leq \overline \vt }   |\vt-\vt_s|^2 } + \delta \intO{   \mathds{1}_{ \vr \leq \underline \vr} },
\nonumber 
\end{align}
and
\begin{align} 
	I_{3,5} & =  \intO{ \mathds{1}_{\underline \vt \leq \vt \leq \overline \vt }  \mathds{1}_{  \vr \geq \overline \vr}\Big( \vt_s^4 - \vt^4 \Big) \Big(T_{k}(\vr) - T_{k}(\vr_s)\Big)} \br  
	& \lesssim  \intO{ \mathds{1}_{\underline \vt \leq \vt \leq \overline \vt }  \mathds{1}_{   \vr \geq \overline \vr}  |\vt-\vt_s| (k+1)} \br 
	& \leq    C(\delta)\intO{ \mathds{1}_{\underline \vt \leq \vt \leq \overline \vt }   |\vt-\vt_s|^2 }   + \delta \intO{   \mathds{1}_{ \vr \geq \overline \vr} \vr^\frac 53 }.
\nonumber 
\end{align}
We conclude,
\begin{align} 
	I_{3} & \leq  C(\delta )\intO{ \mathds{1}_{\vt\geq \overline \vt }   \frac{\kappa (\vt)}{\vt^2} |\Grad \vt |^2 } +  C(\delta)\|\vt-\vt_s\|_{W^{1,2}(\Omega)}^2 \br
	&\quad + \delta \intO{   \mathds{1}_{ \underline \vr \leq \vr \leq \overline \vr}   |\vr-\vr_s|^2 } + \delta \intO{   \mathds{1}_{ \vr \leq \underline \vr} }  +\delta \intO{   \mathds{1}_{ \vr \geq \overline \vr} \vr^\frac 53 }. 
\label{I3-f} 	
\end{align}

\end{itemize}

\medskip

Next we consider the integrals $I_i,i=4,5,6,7,8$ in \eqref{c12-1}.
\begin{itemize}

\item {\bf Integral $I_4$.}	First, 
\begin{align}
		I_4 &= \intO{ \Big(\mathbb{S} (\vt, \Grad \vu)  - \mathbb{S} (\vt_{s}, \Grad \vu_{s})  \Big): \Grad  \phiB} \br
		& = \intO{ \Big(\mathbb{S} (\vt, \Grad \vu)  - \mathbb{S} (\vt, \Grad \vu_s)  + \mathbb{S} (\vt, \Grad \vu_s)  - \mathbb{S} (\vt_{s}, \Grad \vu_{s})  \Big): \Grad  \phiB} \br
		&  \lesssim   \intO{ (1+\vt) | \Grad \vu  - \Grad \vu_s  | |\Grad \phiB | }  +  \intO{ | \Grad \vu_s| |\vt - \vt_s | |\Grad \phiB | }  \br
		& \leq   C(\delta) \intO{  \mathds{1}_{\vt \leq \overline \vt} |\Grad \vu - \Grad \vu_s|^{2} } + \delta \intO{  \mathds{1}_{\vt \leq \overline \vt} |\Grad \phiB|^{2} } \br 
		&\quad  + C(\delta) \intO{  \mathds{1}_{\vt \geq \overline \vt} |\Grad \vu - \Grad \vu_s|^{2} } + C(\delta) \intO{  \mathds{1}_{\vt \geq \overline \vt} \vt^{6} }   + \delta \intO{  \mathds{1}_{\vt \geq \overline \vt} |\Grad \phiB|^{3} } \br
		& \quad + \ep \intO{ |\vt - \vt_s|^{2} }  + \ep\intO{ |\Grad \phiB|^{2} },
	\nonumber 
		\end{align} 
		where we used  $\|\Grad \vu_s\|_{L^\infty(\Omega)}\lesssim \ep $.

Using the boundedness of Bogovskii's operator specified in \eqref{AS6}, we get
\begin{align}
  \intO{  |\Grad \phiB|^{2} }  & \lesssim   \intO{   | T_{k}(\vr) - T_{k}(\vr_{s}) |^{2} }   \br
&   \lesssim   \intO{   \mathds{1}_{ \underline \vr \leq \vr \leq \overline \vr}  (\vr -\vr_{s})^{2} }  + \intO{  \mathds{1}_{ \vr \leq \underline \vr} }  + \intO{  \mathds{1}_{ \vr \geq \overline \vr} (k+1)^{2}}   \br 
& \lesssim \intO{ (\vr - \vr_s)^2 \mathds{1}_{ \underline \vr \leq \vr \leq \overline \vr} } + \intO{ \mathds{1}_{ \vr \leq \underline \vr}  } + \intO{ k \mathds{1}_{  \vr \geq \overline \vr}   \vr^{\frac 53}  } .
\nonumber 
\end{align}
Similarly,
\begin{align}
 \intO{   |\Grad \phiB |^{3} }  & \lesssim   \intO{  | T_{k}(\vr) - T_{k}(\vr_{s}) |^{3} }  \lesssim   \intO{  | T_{k}(\vr) - T_{k}(\vr_{s}) |^{2} }   \br
 &  \lesssim \intO{ (\vr - \vr_s)^2 \mathds{1}_{ \underline \vr \leq \vr \leq \overline \vr} } + \intO{ \mathds{1}_{ \vr \leq \underline \vr}  } + \intO{ k \mathds{1}_{  \vr \geq \overline \vr}   \vr^{\frac 53}  } .  
\nonumber 
\end{align}

Summing the estimates in \eqref{I1-1}--\eqref{I1-4}, together with \eqref{c13}, we conclude
\begin{align}
I_{4} & \leq C(\delta)  \intO{  |\Grad \vu - \Grad \vu_s |^{2} } + C(\delta)  \intO{ \mathds{1}_{\vt \geq \Ov{\vt} } \frac{\kappa (\vt)}{\vt^2} |\Grad \vt |^2 } + C(\delta)  \intO{ |\vt - \vt_s|^{2} }\br
& \quad + \delta \intO{   \mathds{1}_{ \underline \vr \leq \vr \leq \overline \vr}   |\vr-\vr_s|^2 } + \delta \intO{   \mathds{1}_{ \vr \leq \underline \vr} }  +\delta \intO{   \mathds{1}_{ \vr \geq \overline \vr} \vr^\frac 53 }.
\label{I4-f}
\end{align}

\item {\bf Integral $I_5$.} A direct calculation yields
\begin{align} 
I_5   &= \intO{ \Big(\vr \vu \otimes \vu - \vr_{s} \vu_{s}\otimes \vu_s \Big) : \phiB}\br
& = \intO{ \Big(\vr (\vu-\vu_s) \otimes (\vu-\vu_s) + 2 \vr \vu\otimes \vu_s  - \vr \vu_s \otimes \vu_s - \vr_{s} \vu_{s}\otimes \vu_s \Big) : \Grad \phiB} \br 
& = \intO{ \Big(\vr (\vu-\vu_s) \otimes (\vu-\vu_s) + 2 (\vr - \vr_s) (\vu- \vu_s) \otimes \vu_s +  2 \vr_s (\vu- \vu_s) \otimes \vu_s
	\Big) : \Grad \phiB} \br
	 &+ \intO{ ( \vr - \vr_s) \vu_s \otimes \vu_s  \Big) : \Grad \phiB} .
 \nonumber 
\end{align}
First, it is easy to see that  
\begin{align} 
I_{5,1} &  = \intO{ \vr (\vu-\vu_s) \otimes (\vu-\vu_s) 
 : \Grad \phiB }\br
 &  \lesssim \| \vr \|_{L^{\frac{5}{3}}(\Omega)} \| \vu  - \vu_s \|^2_{L^6(\Omega)} \|  T_k(\vr) - T_k(\vr_s) \|_{L^{15}(\Omega)}\lesssim  \|\vu - \vu_s \|^2_{W^{1,2}(\Omega)}.
 \nonumber 
\end{align}

Next,
\begin{align} 
I_{5,2} &  = \intO{  2 (\vr - \vr_s) (\vu- \vu_s) \otimes \vu_s  : \Grad \phiB}\br
&  = \intO{ (\mathds{1}_{\vr\leq \overline \vr} + \mathds{1}_{\vr > \overline \vr} ) 2 (\vr - \vr_s) (\vu- \vu_s) \otimes \vu_s  : \Grad \phiB}\br
 &  \lesssim  \|\vu_s \|_{L^\infty(\Omega)}\| \vu  - \vu_s \|_{L^2(\Omega)} \|  \Grad \phiB \|_{L^{2}(\Omega)} + \|\vu_s \|_{L^\infty(\Omega)} \|\vr  \mathds{1}_{\vr > \overline \vr} \|_{L^{\frac 53}(\Omega)} \| \vu  - \vu_s \|_{L^6(\Omega)} \|  \Grad \phiB \|_{L^{\frac{30}{7}}(\Omega)} \br 
 & \lesssim  \ep  \| \vu -  \vu_s \|^2_{W^{1,2}(\Omega)} + \ep \|  \phiB \|_{L^{2}(\Omega)}^2 + \ep \|\vr  \mathds{1}_{\vr > \overline \vr} \|_{L^{\frac 53}(\Omega)}^2\|  \Grad \phiB \|_{L^{\frac{30}{7}}(\Omega)}^2\br 
 & \lesssim \ep  \| \vu -  \vu_s \|^2_{W^{1,2}(\Omega)} + \ep \intO{   \mathds{1}_{ \underline \vr \leq \vr \leq \overline \vr}   |\vr-\vr_s|^2 } + \ep \intO{   \mathds{1}_{ \vr \leq \underline \vr} }  +\ep\intO{   \mathds{1}_{ \vr \geq \overline \vr} \vr^\frac 53 },
\nonumber 
\end{align}
where we used $\|\vu_s\|_{L^\infty(\Omega)} \lesssim \ep $.

We then calculate
\begin{align} 
I_{5,3} &  = \intO{  2 \vr_s (\vu- \vu_s) \otimes \vu_s  : \Grad \phiB}\br
 &  \lesssim  \|\vu_s \|_{L^\infty(\Omega)}\| \vu  - \vu_s \|_{L^2(\Omega)} \|  \Grad \phiB \|_{L^{2}(\Omega)} \br
 & \lesssim \ep  \|\vu - \vu_s \|^2_{L^2(\Omega)} + \ep \intO{   \mathds{1}_{ \underline \vr \leq \vr \leq \overline \vr}   |\vr-\vr_s|^2 } + \ep \intO{   \mathds{1}_{ \vr \leq \underline \vr} }  +\ep\intO{   \mathds{1}_{ \vr \geq \overline \vr} \vr^\frac 53 }.
\nonumber 
\end{align}

Finally,
\begin{align} 
I_{5,4} &  = \intO{  (\vr - \vr_s)  \vu_s \otimes \vu_s  : \Grad \phiB }\br
 &  \lesssim  \|\vu_s \|_{L^\infty(\Omega)}^2 \Big(\| \mathds{1}_{ \vr \leq \underline \vr} \|_{L^2(\Omega)} \|  \Grad \phiB  \|_{L^{2}(\Omega)}  + \| \mathds{1}_{\underline \vr\leq \vr \leq \overline \vr} (\vr -\vr_s)\|_{L^2(\Omega)} \|  \Grad \phiB \|_{L^{2}(\Omega)} \Big) \br
 &   \quad +  \|\vu_s \|_{L^\infty(\Omega)}^2 \Big(\| \mathds{1}_{\vr \geq \overline \vr} \vr \|_{L^{\frac 53}(\Omega)} \|  \Grad \phiB  \|_{L^{\frac 52}(\Omega)}\Big) \br
 & \lesssim  \ep^2 \intO{   \mathds{1}_{ \underline \vr \leq \vr \leq \overline \vr}   |\vr-\vr_s|^2 } + \ep^2 \intO{   \mathds{1}_{ \vr \leq \underline \vr} }  +\ep^2 \intO{   \mathds{1}_{ \vr \geq \overline \vr} \vr^\frac 53 }.
 \label{I5-4}
\end{align}

We may therefore conclude that 
\begin{align} 
I_{5}   \lesssim  \| \vu -  \vu_s \|^2_{W^{1,2}(\Omega)} + \ep \intO{   \mathds{1}_{ \underline \vr \leq \vr \leq \overline \vr}   |\vr-\vr_s|^2 } + \ep \intO{   \mathds{1}_{ \vr \leq \underline \vr} }  +\ep\intO{   \mathds{1}_{ \vr \geq \overline \vr} \vr^\frac 53 }.
 \label{I5-f}
\end{align}

\item {\bf Integral $I_6$.} Similarly to \eqref{I5-4}, Poincar\'e inequality applied to $\phiB$ yields
 \begin{align} 
I_{6}  & =  - \intO{ (\vr-\vr_{s})\Grad G \cdot \phiB } \br 
&\leq \|\Grad G\|_{L^\infty(\Omega)}\Big( \| \mathds{1}_{ \vr \leq \underline \vr} \|_{L^2(\Omega)} \|   \phiB  \|_{L^{2}(\Omega)}  + \| \mathds{1}_{\underline \vr\leq \vr \leq \overline \vr} (\vr -\vr_s)\|_{L^2(\Omega)} \|  \phiB  \|_{L^{2}(\Omega)} \Big) \br
&  \quad + \|\Grad G\|_{L^\infty(\Omega)}\Big( \| \mathds{1}_{\vr \geq \overline \vr} \vr \|_{L^{\frac 53}(\Omega)} \|  \phiB  \|_{L^{\frac 52}(\Omega)}\Big) \br 
 & \lesssim  \ep \intO{   \mathds{1}_{ \underline \vr \leq \vr \leq \overline \vr}   |\vr-\vr_s|^2 } + \ep \intO{   \mathds{1}_{ \vr \leq \underline \vr} }  +\ep \intO{   \mathds{1}_{ \vr \geq \overline \vr} \vr^\frac 53 }.
 \label{I6-0}
\end{align}

\item {\bf Integral $I_7$.} The renormalized continuity equation \eqref {w2} for the weak solution $(\vr, \vu)$ and for the stationary solution $(\vr_s, \vu_s)$ read
\begin{align}
\partial_{t} T_{k}(\vr) + \Div  ( T_{k}(\vr) \vu)  + (T_{k}'(\vr) \vr - T_{k}(\vr)) \Div \vu = 0,\br  \Div  ( T_{k}(\vr_s) \vu_s)  + (T_{k}'(\vr_s) \vr_s - T_{k}(\vr_s)) \Div \vu_s = 0,
\nonumber 
\end{align}
respectively. Accordingly, 
we may write 
\begin{align} 
I_7  : & = \intO{(\vr \vu - \vr_s \vu_s) \cdot \partial_{t} \phiB }  = - \intO{ (\vr \vu - \vr_s \vu_s) \cdot  \mathcal{B} [ \Div  ( T_{k}(\vr) \vu)  - \Div  ( T_{k}(\vr_s) \vu_s) ]}\br
 & \quad - \intO{ (\vr \vu - \vr_s \vu_s) \cdot \mathcal{B} [\widetilde T_k(\vr, \vu) -\widetilde T_k ( \vr_s,  \vu_s)] } .
\nonumber 
\end{align}
where
\begin{align}
\widetilde T_k(\vr, \vu)  : =  (T_{k}'(\vr) \vr - T_{k}(\vr)) \Div \vu - \langle (T_{k}'(\vr) \vr - T_{k}(\vr)) \Div \vu \rangle.  \nonumber %
\end{align}
Next, we decompose
$$
(\vr \vu - \vr_s \vu_s) =  \vr (\vu - \vu_s) + \vu_s (\vr - \vr_s) , \quad  T_{k}(\vr) \vu  -  T_{k}(\vr_s) \vu_s = T_{k}(\vr) (\vu  - \vu_s) + \vu_s ( T_{k}(\vr) -  T_{k}(\vr_s)) .
 $$

Using the $L^{p}\to L^{p}$ boundedness of the operator $\mathcal{B} \Div$ stated in \eqref{AS7}, we get
\begin{align} 
I_{7,1} & =  -  \intO{ \vr (\vu - \vu_s) \cdot  \mathcal{B} [ \Div  ( T_{k}(\vr) (\vu-\vu_s)) ] } \br 
& \lesssim  \|\vr\|_{L^{\frac 53}(\Omega)} \|\vu - \vu_s\|_{L^{6}(\Omega)} \|T_{k}(\vr) (\vu- \vu_s)\|_{L^{\frac{30}{7}}(\Omega)} \br 
 & \lesssim  \|\vr\|_{L^{\frac 53}(\Omega)} \|\vu - \vu_s\|_{L^{6}(\Omega)} \|\vu - \vu_s\|_{L^{6}(\Omega)}   \|T_{k}(\vr) \|_{L^{15}(\Omega)} \br 
 & \lesssim \| \vu - \vu_s\|^2_{W^{1,2}(\Omega)},
 \nonumber 
\end{align}

\begin{align} 
I_{7,2} & =  -  \intO{ \vr (\vu - \vu_s) \cdot  \mathcal{B} [ \Div  (\vu_s ( T_{k}(\vr) - T_{k}(\vr_s))) ] } \br 
&  \lesssim  \|\vu_s \|_{L^\infty(\Omega)} \|\mathds{1}_{\vr \leq \overline \vr} \vr  \|_{L^\infty(\Omega)}  \| \vu  - \vu_s \|_{L^2(\Omega)} \|  T_{k}(\vr) -  T_{k}(\vr_s) \|_{L^{2}(\Omega)} \br
& \quad  + \|\vu_s \|_{L^\infty(\Omega)} \|\vr  \mathds{1}_{\vr > \overline \vr} \|_{L^{\frac 53}(\Omega)} \| \vu  - \vu_s \|_{L^6(\Omega)} \|  T_{k}(\vr) -  T_{k}(\vr_s)\|_{L^{\frac{30}{7}}(\Omega)} \br 
 & \lesssim  \ep  \| \vu -  \vu_s \|^2_{W^{1,2}(\Omega)} + \ep \|  T_{k}(\vr) -  T_{k}(\vr_s)\|_{L^{2}(\Omega)}^2 + \ep \|\vr  \mathds{1}_{\vr > \overline \vr} \|_{L^{\frac 53}(\Omega)}^2\|  T_{k}(\vr) -  T_{k}(\vr_s) \|_{L^{\frac{30}{7}}(\Omega)}^2\br 
 & \lesssim \ep  \| \vu -  \vu_s \|^2_{W^{1,2}(\Omega)} + \ep \intO{   \mathds{1}_{ \underline \vr \leq \vr \leq \overline \vr}   |\vr-\vr_s|^2 } + \ep \intO{   \mathds{1}_{ \vr \leq \underline \vr} }  +\ep\intO{   \mathds{1}_{ \vr \geq \overline \vr} \vr^\frac 53 },
 \nonumber 
\end{align}

\begin{align} 
I_{7,3} & =  -  \intO{ \vu_s (\vr - \vr_s) \cdot  \mathcal{B} [ \Div  (T_{k}(\vr) (\vu - \vu_s)) ] } \br 
&  \lesssim  \|\vu_s \|_{L^\infty(\Omega)} \|\mathds{1}_{\vr \leq \underline \vr} \|_{L^2(\Omega)}  \| \vu  - \vu_s \|_{L^2(\Omega)} \|  T_{k}(\vr)\|_{L^{\infty}(\Omega)} \br
& \quad  + \|\vu_s \|_{L^\infty(\Omega)} \|\mathds{1}_{\underline \leq \vr \leq \overline \vr} (\vr - \vr_s) \|_{L^2(\Omega)}  \| \vu  - \vu_s \|_{L^2(\Omega)} \|  T_{k}(\vr)\|_{L^{\infty}(\Omega)} \br
& \quad  + \|\vu_s \|_{L^\infty(\Omega)} \|\vr  \mathds{1}_{\vr > \overline \vr} \|_{L^{\frac 53}(\Omega)} \| \vu  - \vu_s \|_{L^6(\Omega)} \|  T_{k}(\vr)\|_{L^{\infty}(\Omega)}\br 
 & \lesssim \ep  \| \vu -  \vu_s \|^2_{W^{1,2}(\Omega)} + \ep \intO{   \mathds{1}_{ \underline \vr \leq \vr \leq \overline \vr}   |\vr-\vr_s|^2 } + \ep \intO{   \mathds{1}_{ \vr \leq \underline \vr} }  +\ep\intO{   \mathds{1}_{ \vr \geq \overline \vr} \vr^\frac 53 },
\nonumber 
\end{align}
and
\begin{align} 
I_{7,4} & =  -  \intO{ \vu_s (\vr - \vr_s) \cdot  \mathcal{B} [ \Div  (\vu_s ( T_{k}(\vr) - T_{k}(\vr_s))) ] } \br 
&  \lesssim  \|\vu_s \|_{L^\infty(\Omega)}^2 \|\mathds{1}_{\vr \leq \underline \vr}  \|_{L^2(\Omega)}   \|  T_{k}(\vr) -  T_{k}(\vr_s) \|_{L^{2}(\Omega)} \br
&  \quad +  \|\vu_s \|_{L^\infty(\Omega)}^2 \|\mathds{1}_{\underline \leq \vr \leq \underline \vr} (\vr -\vr_s) \|_{L^2(\Omega)}    \|  T_{k}(\vr) -  T_{k}(\vr_s) \|_{L^{2}(\Omega)} \br
& \quad  + \|\vu_s \|_{L^\infty(\Omega)}^2 \|\vr  \mathds{1}_{\vr > \overline \vr} \|_{L^{\frac 53}(\Omega)}   \|  T_{k}(\vr) -  T_{k}(\vr_s)\|_{L^{\frac{5}{2}}(\Omega)} \br 
 & \lesssim \ep^2 \intO{   \mathds{1}_{ \underline \vr \leq \vr \leq \overline \vr}   |\vr-\vr_s|^2 } + \ep^2\intO{   \mathds{1}_{ \vr \leq \underline \vr} }  +\ep^2 \intO{   \mathds{1}_{ \vr \geq \overline \vr} \vr^\frac 53 }. 
\nonumber 
\end{align}

\medskip

Similarly, we can handle the remaining integral in $I_7$ obtaining 
\begin{align} 
I_{7,5}  & =  - \intO{ (\vr \vu - \vr_s \vu_s) \cdot \mathcal{B} [\widetilde T_k(\vr, \vu) -\widetilde T_k ( \vr_s,  \vu_s)] } \br 
& \lesssim  \| \vu - \vu_s\|^2_{W^{1,2}(\Omega)} + \ep \intO{   \mathds{1}_{ \underline \vr \leq \vr \leq \overline \vr}   |\vr-\vr_s|^2 } + \ep\intO{   \mathds{1}_{ \vr \leq \underline \vr} }  +\ep \intO{   \mathds{1}_{ \vr \geq \overline \vr} \vr^\frac 53 }.
\nonumber 
\end{align}

We conclude
\begin{align} 
I_{7}  \lesssim  \| \vu - \vu_s\|^2_{W^{1,2}(\Omega)} + \ep \intO{   \mathds{1}_{ \underline \vr \leq \vr \leq \overline \vr}   |\vr-\vr_s|^2 } + \ep\intO{   \mathds{1}_{ \vr \leq \underline \vr} }  +\ep \intO{   \mathds{1}_{ \vr \geq \overline \vr} \vr^\frac 53 }. 
 \label{I7-f}
\end{align}
It is worth noting that we have repeatedly used the uniform boundedness of the norm $\| \vr \|_{L^{\frac{5}{3}}(\Omega)}$ 
stated in \eqref{A6}.

\item {\bf Integral $I_8$.} Finally, we have to check that the integral under the time derivative
\[
I_8 = \intO{ \big(\vr \vu -\vr_{s} \vu_{s} \big)\cdot \phiB }
\]
can be controlled by the relative energy $\intO{ E(\vr, \vt, \vu |\vr_s, \vt_s, \vu_s)}$. To see this, we decompose 
$$
\vr = \vr \mathds{1}_{\vr \leq \underline \vr} + \vr \mathds{1}_{\underline \vr \leq \vr \leq \overline \vr}  + \vr \mathds{1}_{  \vr \geq \overline \vr}, \  \big(\vr \vu -\vr_{s} \vu_{s} \big) = \vr (\vu - \vu_s) + \vu_s (\vr - \vr_s).
$$
Then
\begin{align} 
I_{8,1} & =     \intO{ \vr (\vu - \vu_s) \cdot  \phiB} \br 
&  \lesssim    \|\mathds{1}_{\vr \leq \overline \vr} \vr (\vu - \vu_s)  \|_{L^2(\Omega)}   \| \phiB \|_{L^{2}(\Omega)} +  \|\mathds{1}_{\vr \geq \overline \vr} \vr^{\frac 12} (\vu - \vu_s)  \|_{L^2(\Omega)}  \|\mathds{1}_{\vr \geq \overline \vr} \vr^{\frac 12}   \|_{L^3(\Omega)}   \| \phiB \|_{L^{6}(\Omega)} \br 
 & \lesssim  \| \vr^\frac{1}{2} (\vu - \vu_s)  \|_{L^2(\Omega)}^2 + \|T_k(\vr) - T_k(\vr_s)\|_{L^2(\Omega)}^2\br
 & \lesssim \intO{   \vr  |\vu - \vu_s|^2 }  + \intO{   \mathds{1}_{ \underline \vr \leq \vr \leq \overline \vr}   |\vr-\vr_s|^2 } + \intO{   \mathds{1}_{ \vr \leq \underline \vr} }  + \intO{   \mathds{1}_{ \vr \geq \overline \vr} \vr^\frac 53 },
  \label{I8-1}
\end{align}
and 
\begin{align} 
I_{8,2} & =     \intO{ \vu_s (\vr - \vr_s) \cdot  \phiB} \br 
&  \lesssim    \|  \vu_s\|_{L^{\infty}(\Omega)}  ( \|\mathds{1}_{\vr \leq \underline \vr}  \|_{L^2(\Omega)}  + \|\mathds{1}_{ \underline \vr \leq \vr \leq \overline \vr} (\vr - \vr_s)  \|_{L^2(\Omega)} )  \| \phiB \|_{L^{2}(\Omega)}  +  \|  \vu_s\|_{L^{\infty}(\Omega)} \|\mathds{1}_{\vr \geq \overline \vr} \vr \|_{L^{\frac 53}(\Omega)}   \| \phiB \|_{L^{6}(\Omega)} \br 
 & \lesssim    \ep\Big(\intO{   \mathds{1}_{ \underline \vr \leq \vr \leq \overline \vr}   |\vr-\vr_s|^2 } + \intO{   \mathds{1}_{ \vr \leq \underline \vr} }  + \intO{   \mathds{1}_{ \vr \geq \overline \vr} \vr^\frac 53 }\Big),
  \label{I8-2}
\end{align}
where the integrals on the right-hand sides of \eqref{I8-1} and \eqref{I8-2} can be controlled by the relative energy $\intO{ E(\vr, \vt, \vu |\vr_s, \vt_s, \vu_s)}$.

\end{itemize}

\subsection{Estimates of the remainders}\label{ef}

Our ultimate goal is to control the remainders on the right-hand side of \eqref{c5}.  The first one is
\begin{align} 
  -  \intO{ \Big( \vr (\vu - \vu_s) \otimes (\vu - \vu_s)    \Big) : \Grad \vu_s  } \leq \|\nabla \vu_{s}\|_{L^{\infty}(\Omega)} \|\vu - \vu_{s}\|_{L^{6}(\Omega)}^{2} \|\vr \|_{L^{\frac{3}{2}}(\Omega)} \lesssim  \ep \|\vu - \vu_{s}\|_{W^{1,2}(\Omega)}^{2}. 
\nonumber %
\end{align}

For the second one reads
\begin{align} 
&- \intO{ \frac{\vr-\vr_s}{\vr_s}\Div\mathbb{S}(\vt_s,\Grad\vu_s)\cdot (\vu - \vu_s) } \leq\| \Div\mathbb{S}(\vt_s,\Grad\vu_s) \|_{L^{\infty}(\Omega)}   \intO{ \frac{|\vr-\vr_s|}{\vr_s}  \ |\vu - \vu_s| }  \br 
& \quad \lesssim \ep\Big(\intO{   \mathds{1}_{ \underline \vr \leq \vr \leq \overline \vr}   |\vr-\vr_s|^2 } + \intO{   \mathds{1}_{ \vr \leq \underline \vr} }  + \intO{   \mathds{1}_{ \vr \geq \overline \vr} \vr^\frac 53 } +\|\vu - \vu_{s}\|_{W^{1,2}(\Omega)}^{2} \Big).
\nonumber %
\end{align}

\medskip

Next, we consider the entropy forcing term 
\[
\intO{ \vr (s(\vr, \vt) - s(\vr_s, \vt_s)) (\vu - \vu_{s}) \cdot \Grad \vt_s }.
\]
Seeing \eqref{S3}, we distinguish the radiation component
\[
s_r = \frac{4a}{3 \vr} \vt^3,  
\]
and the ``molecular'' component
\[
s_m(\vr, \vt) = \mathcal{S} \left( \frac{\vr}{\vt^{\frac{3}{2}}} \right). 
\]

The radiation part can be written as
\begin{align}
&\intO{ \vr \left( \frac{\vt^3}{\vr} - \frac{\vt_s^3}{\vr_s} \right) (\vu - \vu_{s})  \cdot \Grad \vt_s } \br
&=  \intO{ \left( \vt^3 - \vt_s^3 \right) (\vu - \vu_{s})  \cdot \Grad \vt_s } + \intO{ \vt_s^3 \frac{ \vr_s - \vr}{\vr_s} (\vu - \vu_{s})  \cdot \Grad \vt_s }\br 
& \lesssim   \|\Grad \vt_s\|_{L^{\infty}(\Omega)}  \Big(\intO{ \mathds{1}_{\vt \leq \overline \vt} (\vt - \vt_s)^{2}   + (\vu - \vu_{s})^{2}} + \intO{ \mathds{1}_{\vt \geq \overline \vt} \vt^{6}    + (\vu - \vu_{s})^{2}} \Big)  \br 
& \quad +\|\Grad \vt_s\|_{L^{\infty}(\Omega)}  \Big( \intO{\mathds{1}_{  \vr \leq \underline \vr}  + \mathds{1}_{ \underline \vr \leq \vr \leq \overline \vr}  ( \vr_s - \vr)^{2}  + (\vu - \vu)^{2} }  + \|  \mathds{1}_{\vr \geq \Ov{\vr} } \vr \|_{L^{\frac{6}{5}}(\Omega)}^2 + \|\vu -\vu_{s}\|_{L^{6}(\Omega)}^{2}\Big)\br 
& \lesssim  \ep\Big( \| \vt - \vt_{s}\|_{W^{1,2}(\Omega)}^{2}  + \intO{   \mathds{1}_{ \vt \geq \overline \vt}  \frac{\kappa(\vt)}{\vt^{2}} |\Grad \vt|^{2}} + \| \vu - \vu_{s}\|_{W^{1,2}(\Omega)}^{2}  \br 
& \quad + \intO{   \mathds{1}_{ \underline \vr \leq \vr \leq \overline \vr}   |\vr-\vr_s|^2 } + \intO{   \mathds{1}_{ \vr \leq \underline \vr} }  + \intO{   \mathds{1}_{ \vr \geq \overline \vr} \vr^\frac 53 }\Big),
\nonumber 
\end{align}
where we have again the fact the weak solution belongs to the bounded absorbing set, in particular,
\[
\|  \mathds{1}_{\vr \geq \Ov{\vr} } \vr \|_{L^{\frac{6}{5}}(\Omega)}^2 \lesssim \|  \mathds{1}_{\vr \geq \Ov{\vr} } \vr \|_{L^{\frac{5}{3}}(\Omega)}^2 \lesssim \|  \mathds{1}_{\vr \geq \Ov{\vr} } \vr \|_{L^{\frac{5}{3}}(\Omega)}^{\frac 53}.
\]

The molecular part reads 
\begin{align}
&\intO{ \vr (s_m(\vr, \vt) - s_m(\vr_s, \vt_s)) (\vu - \vu_{s}) \cdot \Grad \vt_s } \br 
&=  \intO{ \vr (s_m(\vr, \vt) - s_m(\vr_s, \vt)) (\vu - \vu_{s})  \cdot \Grad \vt_s } + 
\intO{ \vr (s_m(\vr_s, \vt) - s_m(\vr_s, \vt_s)) (\vu - \vu_{s})  \cdot \Grad \vt_s }. 
\nonumber
\end{align}
It follows the Third law of thermodynamics enforced by hypothesis \eqref{S7} that
\[
\vr |s_m (\vr, \vt)| \lesssim 1 + \vr \log^+(\vr) + \vr \log^+(\vt).
\]
Consequently, the desired bounds follow from the existing estimates \eqref{A6} on bounded absorbing set:
 \begin{align}
&\intO{ \vr (s_m(\vr, \vt) - s_m(\vr_s, \vt_s)) (\vu - \vu_{s}) \cdot \Grad \vt_s } \br 
&\lesssim  \ep\Big( \| \vt - \vt_{s}\|_{W^{1,2}(\Omega)}^{2}  + \intO{   \mathds{1}_{ \vt \geq \overline \vt}  \frac{\kappa(\vt)}{\vt^{2}} |\Grad \vt|^{2}} + \| \vu - \vu_{s}\|_{W^{1,2}(\Omega)}^{2}  \br 
& \quad + \intO{   \mathds{1}_{ \underline \vr \leq \vr \leq \overline \vr}   |\vr-\vr_s|^2 } + \intO{   \mathds{1}_{ \vr \leq \underline \vr} }  + \intO{   \mathds{1}_{ \vr \geq \overline \vr} \vr^\frac 53 }\Big).
\nonumber 
\end{align}

\medskip

The rest of the remainder terms can be estimated in a similar manner and we shall not repeat the details. Consequently, the right-hand side of \eqref{c5}  can be controlled by 
 \begin{align}
 &\ep\Big( \| \vt - \vt_{s}\|_{W^{1,2}(\Omega)}^{2}  + \intO{   \mathds{1}_{ \vt \geq \overline \vt}  \frac{\kappa(\vt)}{\vt^{2}} |\Grad \vt|^{2}} + \| \vu - \vu_{s}\|_{W^{1,2}(\Omega)}^{2}  \br 
& \quad + \intO{   \mathds{1}_{ \underline \vr \leq \vr \leq \overline \vr}   |\vr-\vr_s|^2 } + \intO{   \mathds{1}_{ \vr \leq \underline \vr} }  + \intO{   \mathds{1}_{ \vr \geq \overline \vr} \vr^\frac 53 }\Big).
\nonumber
\end{align}

\subsection{Conclusion, unconditional convergence}

As a consequence of the estimates obtained in Section \ref{D}, we may
choose $\delta$ suitably small so that 
\begin{align} 
&\intO{ \Big( p(\vr, \vt_s) - p(\vr_s, \vt_s) \Big) \Big(T_{k}(\vr) - T_{k}(\vr_s)\Big) } \br
 & \lesssim  -  \frac{\D}{\dt} \intO{ \Big(\vr \vu -\vr_{s} \vu_{s} \Big)\cdot \phiB }  +  \Big(\|\vu-\vu_s\|_{W^{1,2}(\Omega)}^2 +   \|\vt-\vt_s\|_{W^{1,2}(\Omega)}^2 + \intO{ \mathds{1}_{\vt\geq \overline \vt }   \frac{\kappa (\vt)}{\vt^2} |\Grad \vt |^2 }   \Big).
\label{e1}
\end{align} 
Multiplying inequality \eqref{e1} by a sufficiently small positive constant and adding the result to \eqref{c5}, using the lower bound estimates on the temperature and velocity dissipations \eqref{cc14}, together with the estimates obtained in Section \ref{ef},  we get
\begin{align}
	&\frac{\D}{\dt} \intO{ \left[  E \left( \vr, \vt, \vu \ \Big| \vr_s , \vt_s, \vu_{s} \right) + \delta \Big(\vr \vu -\vr_{s} \vu_{s} \Big)\cdot \phiB \right]  } 	+ 2  \delta  \| \vu - \vu_{s}\|_{W^{1,2}(\Omega)}^{2} \br
	&\quad + 2 \delta \Big(\| \vt - \vt_s \|^2_{W^{1,2}_0 (\Omega)}  + \intO{   \frac{\kappa(\vt)}{\vt^2}  |\Grad \vt|^2 } +  \intO{  \mathds{1}_{\vt \geq \underline{\vt}}  \frac{\kappa(\vt)}{\vt^2}  |\Grad \vt - \Grad \vt_s|^2 } \Big)\br
	& \quad +  2\delta \Big(  \intO{   \mathds{1}_{ \underline \vr \leq \vr \leq \overline \vr}   |\vr-\vr_s|^2 }  + \intO{   \mathds{1}_{ \vr \leq \underline \vr} }  + \intO{   \mathds{1}_{ \vr \geq \overline \vr} \vr^\frac 53 }\Big)\br
	& \lesssim  \ep\Big( \| \vt - \vt_{s}\|_{W^{1,2}(\Omega)}^{2}  + \intO{   \mathds{1}_{ \vt \geq \overline \vt}  \frac{\kappa(\vt)}{\vt^{2}} |\Grad \vt|^{2}} + \| \vu - \vu_{s}\|_{W^{1,2}(\Omega)}^{2}  \br 
& \quad + \intO{   \mathds{1}_{ \underline \vr \leq \vr \leq \overline \vr}   |\vr-\vr_s|^2 } + \intO{   \mathds{1}_{ \vr \leq \underline \vr} }  + \intO{   \mathds{1}_{ \vr \geq \overline \vr} \vr^\frac 53 }\Big),
	\nonumber 
\end{align}
for some $\delta>0$. Thus if $0<\ep \leq \ep_{0}$ is sufficiently small, there holds
\begin{align}
	&\frac{\D}{\dt} \intO{ \left[  E \left( \vr, \vt, \vu \ \Big| \vr_s , \vt_s, \vu_{s} \right) + \delta \Big(\vr \vu -\vr_{s} \vu_{s} \Big)\cdot \phiB \right]  } 	+ \delta  \| \vu - \vu_{s}\|_{W^{1,2}(\Omega)}^{2}  \br
	&\quad + \delta \Big(\| \vt - \vt_s \|^2_{W^{1,2}_0 (\Omega)}  + \intO{   \frac{\kappa(\vt)}{\vt^2}  |\Grad \vt|^2 } +  \intO{  \mathds{1}_{\vt \geq \underline{\vt}}  \frac{\kappa(\vt)}{\vt^2}  |\Grad \vt - \Grad \vt_s|^2 } \Big)\br
	& \quad +  \delta \Big(  \intO{   \mathds{1}_{ \underline \vr \leq \vr \leq \overline \vr}   |\vr-\vr_s|^2 }  + \intO{   \mathds{1}_{ \vr \leq \underline \vr} }  + \intO{   \mathds{1}_{ \vr \geq \overline \vr} \vr^\frac 53 }\Big)\br
	& \leq 0.
	\label{e3}
\end{align}
Note that the estimates \eqref{I8-1}--\eqref{I8-2} ensure 
\[
\delta \Big(\vr \vu -\vr_{s} \vu_{s} \Big)\cdot \phiB  \leq \frac{1}{2} E \left( \vr, \vt, \vu \ \Big| \vr_s , \vt_s, \vu_{s} \right)
\]
provided  $\delta>0$ is suitably small. Therefore, our desired unconditional stability result \eqref{M4} follows from \eqref{e3}.
We have proved Theorem \ref{TM1}.

\section{Concluding discussion, the Rayleigh--B\' enard problem}
\label{cr}

The proof of Theorem \ref{TM1} can be easily modified to establish unconditional convergence to equilibrium for the (compressible) \emph{Rayleigh-B\' enard} problem. 
Specifically, the motion of a fluid is confined between
two parallel planes, where the temperature of the bottom plane is maintained at the
level $\Theta_B$, while the top plane has the ambient temperature $\Theta_U$, typically $\Theta_B > \Theta_U$. The fluid is subjected to the gravitational force 
acting in the downward direction. The problem is supplemented by the periodic boundary conditions in the horizontal variables $(x_1,x_2)$. 
Consequently, the relevant boundary conditions read 
\begin{align} 
\Omega &= \mathbb{T}^2 \times (0,1),\ \mathbb{T}^2 \equiv \left( [-1,1]|_{\{ -1,1\}} \right)^2, \br
\vu|_{\partial \Omega} &= 0, \br
\vt|_{x_3 = 0} &= \Theta_B,\ \vt|_{x_3 = 1} = \Theta_U.
	\label{cr1}
	\end{align}
The gravitational force acts in the direction $\vc{g}$, meaning 
\begin{equation} \label{cr2}
G = \vc{g} \cdot x.
\end{equation}

If $\Theta_U$, $\Theta_B$ are positive constants, 
it is easy to check that the stationary problem \eqref{p9}--\eqref{p11} admits a static solution ($\vu_s = 0$) only if 
$\vc{g} = (0,0,g)$, meaning the gravitational force acts in the vertical direction, cf. Daniels et al. \cite{DBPB}.  

The proof of Theorem \ref{TM1} can be repeated with only formal modifications to obtain the following result. 
\begin{Theorem}[\bf Rayleigh--B\' enard problem -- unconditional stability] \label{TM2}
	Let 
	\[
	\Omega = \mathbb{T}^2 \times (0,1),\ \mathbb{T}^2 := \left( [-1,1]|_{\{ -1,1\}} \right)^2.
	\]
	Let the boundary conditions by given by \eqref{cr1}, where $\Theta_B$, $\Theta_U$ are positive constants. Let 
	\[
	G = \vc{g} \cdot x,\ \vc{g} \in R^3.
	\]
	 Suppose the thermodynamic 
	functions $p$, $e$, and $s$ satisfy the hypotheses \eqref{S1}--\eqref{convex-pm}, while the transport coefficients 
	$\mu$, $\eta$, and $\kappa$ comply with \eqref{S10}--\eqref{S12}. Let 
	the total mass of the fluid $m_0 > 0$ be given. 
	
	Then there exists $\ep > 0$ such that 
	\[
		\intO{ E \left( \vr, \vt, \vu \Big| \vr_s , \vt_s, \vu_s \right) (t, \cdot) } \to 0 
		\ \mbox{as}\ t \to \infty
	\]
	for any global--in--time weak solution $(\vr, \vt, \vu)$ of the NSF system defined on $(T, \infty)$ with
	\[
	\intO{\vr } = m_0,
	\]
where $(\vr_s , \vt_s, \vu_s)$ is the solution of the stationary problem \eqref{p9}--\eqref{p11} provided
\[
|\Theta_B - \Theta_U| + |\vc{g}| \leq \ep.	
\]	
\end{Theorem}

Once more, we point out that the convergence stated in Theorem \ref{TM2} is global, meaning it holds for any weak solution of the NSF system. Results concerning local stability of stationary states for the Rayleigh--B\' enard problem driven by a vertical gravitational force were obtained by Nishida, Padula, Teramoto \cite{NiPaTeI}, \cite{NiPaTeII}, \cite{NiPaTeIII}.

\def\cprime{$'$} \def\ocirc#1{\ifmmode\setbox0=\hbox{$#1$}\dimen0=\ht0
	\advance\dimen0 by1pt\rlap{\hbox to\wd0{\hss\raise\dimen0
			\hbox{\hskip.2em$\scriptscriptstyle\circ$}\hss}}#1\else {\accent"17 #1}\fi}


\end{document}